\documentclass[11pt]{article}

\usepackage{amsmath}

\usepackage{amssymb}
\usepackage{latexsym}
  \newcommand{\argmax}{\operatornamewithlimits{argmax}}

\begin{document}

\newtheorem{theo}{Theorem}
\newtheorem{prop}{Proposition}
\newtheorem{lemma}{Lemma}
\newtheorem{cor}{Corollary}
\newtheorem{defi}{Definition}
\newtheorem{Rk}{Remark}
\begin{center}
{ \Large \bf \sc Joint behaviour of semirecursive kernel
\vspace{0.25cm}\\ estimators of the location and of the size of the mode 
\vspace{0.25cm}\\
of a probability density function
}

\vspace{1.cm}
Abdelkader Mokkadem \hspace{1.2cm} Mariane Pelletier \hspace{1.2cm} Baba Thiam

\vspace{0.3cm}
(mokkadem, pelletier, thiam)@math.uvsq.fr

\vspace{0.5cm}

{\small Universit\'e de Versailles-Saint-Quentin\\
D\'epartement de Math\'ematiques\\
45, Avenue des Etats-Unis\\
78035 Versailles Cedex\\
France }
\vspace{1.cm}

{\bf Abstract}:\ \parbox[t]{9.cm}{\small
 Let $\theta$ and $\mu$ denote the location and the size of the mode of a 
probability density. We study the joint convergence rates of semirecursive kernel 
estimators of $\theta$ and $\mu$. We show how the estimation of the size of the 
mode allows to measure the relevance of the estimation of its location. 
We also enlighten that, beyond their computational advantage on nonrecursive 
estimators, the semirecursive estimators are preferable to use for the construction 
on confidence regions.} 
\vspace{2.cm}

\end{center}

\noindent
{\bf AMS Subj. Classification:} 62G07, 62G20\\
\noindent
{\bf Key words~:}
Location and size of the mode; semirecursive estimation; central limit theorem; 
law of the iterated logarithm
\vspace{0.5cm}

\newpage
\section{Introduction}

Let $X_{1},\dots,X_{n}$  be  independent and
identically distributed $\mathbb{R}^d$-valued random variables with unknown 
probability density $f$. 
The aim of this paper is to study the joint kernel estimation of the location 
$\theta$ and of the size $\mu=f(\theta)$ of the mode of $f$. The mode is 
assumed to be unique, that is,  $f(x)<f(\theta)$ for any $x \neq \theta$, and 
nondegenerated, that is, the second order differential $D^2f(\theta)$ 
at the point $\theta$ is nonsingular (in the sequel, $D^mg$ will 
denote the differential of order $m$ of a multivariate function $g$). \\

The problem of estimating the location of the mode of a probability density was
widely studied. Kernel methods were considered, among many others, by 
Parzen \cite{parzenm}, Nadaraya \cite{nadarayam}, Van Ryzin \cite{vanryzinm}, 
R\"uschendorf \cite{ruschendorfm}, Konakov \cite{konakovm}, 
Samanta  \cite{samantam}, Eddy (\cite{eddy1m}, \cite{eddy2m}), 
Romano \cite{romanom}, Tsybakov \cite{tsybakov}, Vieu \cite{vieum}, 
Mokkadem and Pelletier \cite{mokkadempelletierm}, and 
Abraham et al. (\cite{abraham1}, \cite{abraham2}). 
At our knowledge, the behaviour of estimators of the size of the mode has not 
been investigated in detail, whereas there are at least two statistical 
motivations for estimating this parameter. First, a use of an estimator 
of the size is necessary for the construction of confidence 
regions for the location of the mode (see, e.g., Romano \cite{romanom}). 
As a more important motivation, let us underline that the high of the peak 
gives information on the shape of a density; from this point view, as 
suggested by Vieu \cite{vieum}, the location of the mode is more related 
to the shape of the derivative of $f$, whereas the size of the mode is 
more related to the shape of the density itself. Moreover, the knowledge of 
the size of the mode allows to measure the pertinence of the parameter 
location of the mode.
 \\

Let us mention that, even if the problem of estimating the size of the mode 
was not investigated in the framework of density estimation, it was studied 
in the framework of regression estimation. M\"uller \cite{muller} proves in particular the joint 
asymptotic normality and independence of kernel estimators of the location and 
of the size of the mode in the framework of nonparametric regression models 
with fixed design. In the framework of nonparametric regression 
with random design, a similar result is obtained by Ziegler (\cite{ziegler a}, \cite{ziegler b})
for kernel estimators, 
and by Mokkadem and Pelletier \cite{mokka} for estimators issued from 
stochastic approximation methods. \\

This paper is focused on semirecursive kernel estimators of $\theta$ and 
$f(\theta)$. To explain why we chose this option of semirecursive estimators, 
let us first recall that the (nonrecursive) wellknown kernel estimator of the location 
of the mode introduced by Parzen \cite{parzenm} is defined 
as a random variable $\theta^*_n$ satisfying 
\begin{eqnarray*}
f_n^*(\theta_n^*)=\sup_{y \in \mathbb{R}^d}f_n^*(y),
\end{eqnarray*}
where $f_n^*$ is Rosenblatt's estimator of $f$; more precisely, 
\begin{eqnarray*}
f_n^*(x)=\frac{1}{nh_n^d}\sum_{i=1}^{n}K\left(\frac{x-X_i}{h_n}\right),
\end{eqnarray*}
where the bandwidth $(h_n)$ is a sequence of positive real numbers going to zero and 
the kernel $K$ is a continuous function satisfying $\lim_{\|x\| \to +\infty}K(x)=0$, 
$\int_{\mathbb{R}^d} K(x)dx=1$. 
The asymptotic behaviour of $\theta_n^*$ was widely studied (see, among others, 
\cite{eddy1m}, \cite{eddy2m}, 
\cite{konakovm}, 
\cite{mokkadempelletierm}, 
\cite{nadarayam}, 
\cite{parzenm},
\cite{romanom},
\cite{ruschendorfm},  
\cite{samantam}, 
\cite{vanryzinm}, 
\cite{vieum}), but, on a computational 
point of view, the estimator $\theta^*_n$ has a main drawback: its update, 
from a sample size $n$ to a sample size $n+1$, is far from being immediate. 
Applying the stochastic approximation method, 
Tsybakov \cite{tsybakov} introduced the recursive kernel estimator of 
$\theta$ defined as 
$$T_n=T_{n-1}+\gamma_n\frac{1}{h_n^{d+1}}
\nabla K\left(\frac{T_{n-1}-X_n}{h_n}\right),$$
where $T_0\in\mathbb{R}^d$ is arbitrarily chosen, and  
the stepsize $(\gamma_n)$ is a sequence of positive real numbers going to zero. 
The great property of this estimator is that its update is very rapid. 
Unfortunately, for reasons inherent to stochastic approximation algorithms 
properties, very strong assumptions on the density $f$ must be required 
to ensure its consistency. 
A recursive version $f_n$ of Rosenblatt's density estimator was introduced 
by Wolverton and Wagner \cite{wolvertonm} (and discussed, among others, 
by Yamato \cite{yamatom}, Davies \cite{daviesm}, Devroye \cite{devroyem}, 
Menon et al. \cite{menonprasadsinghm}, Wertz \cite{wertzm}, 
Wegman and Davies \cite{wegmandaviesm}, Roussas \cite{roussasm}, 
and Mokkadem et al. \cite{mokkadempelletierthiamm}). Let us recall that $f_n$ 
is defined as 
\begin{eqnarray}
\label{estmode}
 f_n(x)=\frac{1}{n}\sum_{i=1}^{n} \frac{1}{h_i^d}K\left(\frac{x-X_i}{h_i}\right).
\end{eqnarray} 
Its update from a sample of size $n$ to one of size $n+1$ is immediate since 
$f_n$ clearly satisfies the recursive relation
$$f_n(x)=\left(1-\frac{1}{n}\right)f_{n-1}(x)+\frac{1}{nh_n^d}K\left(\frac{x-X_n}{h_n}\right).$$
This property of rapid update of the density estimator is particularly important 
in the framework of mode estimation, since the number of points where $f$ must be estimated 
is very large. We thus define a semirecursive version of Parzen's estimator of the location 
of the mode by using Wolverton-Wagner's recursive density estimator, rather than 
Rosenblatt's density estimator. More precisely, our estimator $\theta_n$ of the location 
$\theta$ of the mode is a random variable satisfying 
\begin{eqnarray}
\label{babamode}
 f_n(\theta_n)=\sup_{y \in \mathbb{R}^d}f_n(y).
\end{eqnarray} 

Let us mention that, in the same way as for Parzen's estimator, the fact that the kernel 
$K$ is continuous and vanishing at infinity ensures that the choice of $\theta_n$ 
as a random variable satisfying (\ref{babamode}) can be made with the help 
of an order on $\mathbb{R}^d$. For example, one can consider the following 
lexicographic order: $x \leq y$ if the first nonzero coordinate of $x-y$ is negative. 
The definition
\begin{eqnarray*}
\theta_n = \inf\left \lbrace y \in \mathbb{R}^d \ 
\mbox{such that} \ f_n(y)=\sup_{x \in \mathbb{R}^d} f_n(x)\right\},
\end{eqnarray*}
where the infimum is taken with respect to the lexicographic order on $\mathbb{R}^d$, 
ensures the measurability of the kernel mode estimator. \\

Let us also mention that, in order to make more rapid the computation of the kernel 
estimator of the location of the mode, Abraham et al. (\cite{abraham1}, \cite{abraham2}) 
proposed the following alternative version of Parzen's estimator $\theta_n^*$: 
\begin{eqnarray*}
\label{modebiau}
\hat\theta_n^*=\argmax_{1\leq i\leq n}f_n^*(X_i).
\end{eqnarray*}
Similarly, we could consider the following alternative version of our semirecursive 
estimator $\theta_n$:
$$\hat\theta_n=\argmax_{1\leq i\leq n}f_n(X_i).$$
However, to establish the asymptotic properties of $\hat\theta_n^*$, Abraham et al. 
\cite{abraham2} prove the asymptotic proximity between $\theta_n^*$ and $\hat\theta_n^*$, 
which allows them to deduce the asymptotic weak behaviour of $\hat\theta_n^*$ from the 
one of $\theta_n^*$. In the same way, we can conjecture that the asymptotic weak 
behaviour of $\hat\theta_n$ could be deduced from the one of $\theta_n$, but, in this paper, 
we limit ourselves on establishing the asymptotic properties of $\theta_n$.\\

Let us now come back to the problem of estimating the size $f(\theta)$ of the mode. The ordinarily 
used estimator is defined as $\mu_n^*=f_n^*(\theta_n^*)$ ($f_n^*$ being Rosenblatt's density 
estimator and $\theta_n^*$ Parzen's mode estimator); the consistency of $\mu_n^*$ is sufficient 
to allow the construction of confidence regions for $\theta$ (see, e.g., Romano \cite{romanom}). 
Adapting the construction of $\mu_n^*$ to the semirecursive framework would lead us to estimate 
$f(\theta)$ by 
\begin{eqnarray}
\label{bebemode}
\mu_n=f_n(\theta_n).
\end{eqnarray}
However, this estimator has two main drawbacks (as well as $\mu_n^*$). 
First, the use of a higher order 
kernel $K$ is necessary for $(\mu_n-\mu)$ to satisfy a central limit theorem, and 
thus for the construction of confidence intervals of $\mu$ (and of confidence regions 
for $(\theta,\mu)$). Moreover, in the case when a higher order kernel is used, it 
is not possible to choose a bandwidth for which both estimators $\theta_n$ and 
$\mu_n$ converge at the optimal rate. These constations lead us to use two 
different bandwidths, one for the estimation of $\theta$, the other one for the 
estimation of $\mu$. More precisely, let $\tilde{f}_n$ be the recursive kernel 
density estimator defined as
\begin{eqnarray*}
\tilde{f}_n(x)=\frac{1}{n}\sum_{i=1}^n\frac{1}{\tilde{h}_i^d}
K\left(\frac{x-X_i}{\tilde{h}_i}\right),
\end{eqnarray*} 
where the bandwidth $(\tilde{h}_n)$ may be different from $(h_n)$ used in the 
definition of $f_n$ (see (\ref{estmode})); we estimate the size of the mode by 
\begin{eqnarray}
\label{bis}
\tilde{\mu}_n = \tilde{f}_n(\theta_n),
\end{eqnarray}
where $\theta_n$ is still defined by (\ref{babamode}), and thus with the first 
bandwidth $(h_n)$. \\

The purpose of this paper is the study of the joint asymptotic behaviour of 
$\theta_n$ and $\tilde\mu_n$. We first prove the strong consistency of both 
estimators. We then establish the joint weak convergence rate of 
$\theta_n$ and $\tilde\mu_n$. We prove in particular that 
adequate choices of the bandwidths lead to the asymptotic normality and 
independence of these estimators, and that the use of different bandwidths 
allow to obtain simultaneously the optimal convergence rate of both estimators. 
We then apply our weak convergence rate result to the construction of confidence 
regions for $(\theta,\mu)$, and illustrate this application with a simulations 
study. This application enlightens the advantage of using  semirecursive 
estimators rather than nonrecursive estimators. It also shows how the estimation 
of the size of the mode gives information on the relevance of estimating its 
location. Finally, we establish 
the joint strong  convergence rate of $\theta_n$ and $\tilde\mu_n$.

\section{Assumptions and Main Results}

Throughout this paper, $(h_n)$ and $(\tilde{h}_n)$ 
are defined as $h_n=h(n)$ and $\tilde{h}_n=\tilde{h}(n)$ for all $n\geq 1$, 
where $h$ and $\tilde{h}$ are two positive functions.

\subsection{Strong consistency}
\label{2.1}

The conditions we require for the strong consistency of $\theta_n$ and $\tilde{\mu}_n$ 
are the following.
\begin{description}
\item (A1) i) $K$ is an integrable, differentiable, and even function such that 
$\int_{\mathbb{R}^d}K(z)dz =1$.\\
          ii) There exists $\zeta>0$ such that $\int_{\mathbb{R}^d}\|z\|^{\zeta}|K(z)dz|<\infty$.\\
         iii) $K$ is H\"older continuous.\\                 
          iv) There exists $\gamma>0$ such that $z \mapsto\|z\|^{\gamma}\left|K(z)\right|$ is a bounded function.
\item (A2) i) $f$ is uniformly continuous on $\mathbb{R}^d$.\\
           ii) There exists $\xi>0$ such that $\int_{\mathbb{R}^d}\|x\|^{\xi}f(x)dx<\infty$.\\
          iii) There exists $\eta>0$ such that $z\mapsto\|z\|^{\eta}f(z)$ is a bounded function.\\
           iv) There exists $\theta \in \mathbb{R}^d$ such that $f(x)<f(\theta)$ for all $x \neq \theta$.
\item (A3) The functions $h$ and $\tilde{h}$ are locally bounded and vary regularly with exponent $(-a)$ and $(-\tilde{a})$ respectively, where $\displaystyle a\in ]0,1/(d+4)[$, $\displaystyle\tilde{a}\in \left]0,1/(d+2)\right[$.  
\end{description}
\begin{Rk}\label{alpha}
Note that (A1)iv) implies that $K$ is bounded.
\end{Rk}
\begin{Rk}\label{+delta}
The assumptions required on the probability density to establish the strong 
consistency of the semirecursive estimator of the location of the mode are 
slightly stronger than those needed for the nonrecursive estimator (see, e.g., 
\cite{mokkadempelletierm}, \cite{romanom}), but are much weaker than the ones 
needed for the recursive estimator (see \cite{tsybakov}). 
\end{Rk}
\begin{Rk}\label{beta}
Let us recall that a positive function (not necessarily monotone) $\mathcal{L}$ defined on $]0,\infty[$ is slowly varying if  $\lim_{t \to \infty}\mathcal{L}(tx)/\mathcal{L}(t)=1$, and that a function $G$ varies regularly with exponent $\rho$, $\rho \in \mathbb{R}$, if and only if it is of the form $G(x)=x^{\rho}\mathcal{L}(x)$ with $\mathcal{L}$ slowly varying (see, for example, Feller \cite{fellerm} page 275). Typical examples of regularly varying functions are $x^\rho$, $x^\rho\log x$, $x^\rho\log \log x$, $\displaystyle x^\rho\log x/\log \log x$, and so on.
\end{Rk}

\begin{prop}\label{Prop}
Let $\theta_n$ and $\tilde\mu_n$ be defined by (\ref{babamode}) and (\ref{bis}), 
respectively. Under (A1)-(A3), 
\begin{eqnarray*}
\lim_{n \to \infty}\theta_n=\theta \ \ \mbox{a.s.} \ \ \mbox{and}\ \ \lim_{n \to \infty} \tilde{\mu}_n=\mu \ \ \mbox{a.s.} 
\end{eqnarray*}
\end{prop}

\subsection{Weak convergence rate}
\label{2.2}

In order to state the weak convergence rate of $\theta_n$ and $\tilde{\mu}_n$, 
we need the following additional assumptions on $K$ and $f$.
\begin{description}
\item (A4) i) $K$ is twice differentiable on $\mathbb{R}^d$.\\ 
          ii) $z\mapsto z\nabla K(z)$ is integrable.\\
         iii) For any $(i,j) \in \left\{1,\ldots,d\right\}^2$, $\partial^{2}K/\partial x_i\partial x_j$ is bounded integrable and H\"older continuous.\\
iv) $K$ is a kernel of order $q \geq 2$ i.e. $\forall s \in \{1,\ldots,q-1\}$, $\forall j \in \{1,\ldots,d\}$, $\int_{\mathbb{R}^d}y_j^sK(y)dy_j=0$ and $\int_{\mathbb{R}^d}|y_j^qK(y)|dy<\infty$.\\
 \item (A5) i) $D^2f(\theta)$ is nonsingular.\\
           ii) $D^2f$ is $q$-times differentiable, $\nabla f$ and $D^qf$ are bounded.\\
          iii) For any $(i,j) \in \left\{1,\ldots,d\right\}^2$, $\sup_{x\in \mathbb{R}^d}\|D^q\left(\partial^{2}f/\partial x_i\partial x_j\right)\|<\infty$, and for any $k \in \left\{1,\ldots,d\right\}$, $\sup_{x\in \mathbb{R}^d}\|D^q(\partial f/\partial x_k)\|<\infty$. 
\end{description}
\begin{Rk}\label{delta}
Note that (A4)ii) and (A4)iii) imply that $\nabla K$ is Lipschitz-continuous and integrable;\\ it is thus straightforward to see that $\lim_{\|x\| \to \infty}\|\nabla K(x)\|=0$ (and in particular $\nabla K$ is bounded).
\end{Rk}
We also need to add conditions on the bandwiths. Let us set
\begin{eqnarray*}
\mathcal{L}_{\theta}(n)=n^{a}h_n \ \  \mbox{and}\ \  \mathcal{L}_{\mu}(n)=n^{\tilde a}\tilde h_n.
\end{eqnarray*}
(In view of (A3), $\mathcal{L}_{\theta}$ and $\mathcal{L}_{\mu}$ are positive slowly 
varying functions, see Remark \ref{beta}). In the statement of the 
the weak convergence rate of $\theta_n$ and $\tilde{\mu}_n$, we shall refer to the 
following conditions. 
\begin{description} 
\item (C1) One of the following two conditions is fulfilled. 
\begin{description}
\item   i) $\dfrac{1}{d+4} < \tilde{a} <\dfrac{q}{d+2q+2}$ and $\dfrac{\tilde a}{q}<a<\dfrac{1-2\tilde a}{d+2}$;
\item  ii) $\dfrac{1}{d+2q}<\tilde a \leq\dfrac{1}{d+4}$ and $\dfrac{1}{d+2q+2}<a<\dfrac{1+\tilde ad}{2\left(d+2\right)}$ 
\end{description}
\item (C2)  One of the following two conditions is fulfilled. 
\begin{description}
\item  i) $0< \tilde{a} <  \dfrac{1}{d+2q}$ and  $\dfrac{\tilde a}{2}<a<\dfrac{1}{d+2q+2}$;
\item ii) $\tilde a = \dfrac{1}{d+2q}$, $\lim_{n \to \infty}\mathcal{L}_{\mu}(n)=\infty$ and $\dfrac{1}{2(d+2q)}<a <  \dfrac{1}{d+2q+2}$.
\end{description}
\end{description}
\begin{Rk}
(C1) implies that  
$\lim_{n \to \infty}{nh_n^{d+2q+2}}=0$ and $\lim_{n \to \infty}{n\tilde h_n^{d+2q}}=0$, 
whereas (C2) implies that  
$\lim_{n \to \infty}{nh_n^{d+2q+2}}=\infty$ and $\lim_{n \to \infty}{n\tilde h_n^{d+2q}}=\infty$.
\end{Rk}
We finally need to introduce the following notation:
\begin{eqnarray}
\label{recmode}
B_q(\theta) =  \left(
\begin{array}{c}
\frac{(-1)^{q}}{q!(1-aq)}\nabla\Big(\sum_{j=1}^d\beta_j^q\frac{\partial^qf}{\partial{x_j^q}}(\theta)\Big)\\
\frac{(-1)^q}{q!(1-\tilde aq)}\sum_{j=1}^d\beta_j^q\frac{\partial^qf}{\partial{x_j^q}}(\theta)
\end{array}\right) \quad \mbox{with} \quad \beta_j^q=\int_{\mathbb{R}^d}y_j^qK(y)dy, \ \  aq\neq 1\ \ \mbox{and}\ \  \tilde aq\neq 1,
\end{eqnarray}
\begin{eqnarray}
\label{racmode}
A = \left(
\begin{array}{c c}
-\left[D^2f(\theta)\right]^{-1} & 0 \\
0 & 1
\end{array}\right), \ \  \Sigma = \left(
\begin{array}{cc}
\frac{f(\theta)G}{1+a(d+2)} & 0 \\
0 & \frac{f(\theta)\int_{\mathbb{R}^d}K^2(z)dz}{1+\tilde ad}
\end{array}\right),
\end{eqnarray}
$G$ is the matrix $d \times d$ defined by $\displaystyle G^{(i,j)}= \int_{\mathbb{R}^d}\frac{\partial{K}}{\partial{x_i}}(x)\frac{\partial{K}}{\partial{x_j}}(x)dx$, and, for any 
$c, \tilde{c} \geq 0$, $D(c,\tilde{c})=\left(
\begin{array}{cc}
\sqrt{c}I_d & 0 \\
0 & \sqrt{\tilde{c}}
\end{array}\right)$ where $I_d$ is the $d\times d$ identity matrix.\\

\begin{theo}\label{cen}
Let $\theta_n$ and $\tilde\mu_n$ be defined by (\ref{babamode}) and (\ref{bis}), respectively,
and assume that (A1)-(A5) hold. 
\begin{description}
\item  i) If (C1) 
is satisfied, then
\begin{eqnarray*}
\left(
\begin{array}{c}
\sqrt{nh_n^{d+2}}(\theta_n-\theta)\\
\sqrt{n\tilde{h}_n^{d}}(\tilde{\mu}_n-\mu)
\end{array}\right)\stackrel{\mathcal{D}}{\longrightarrow}                             \mathcal{N} 
\left(0, A\Sigma A\right).
\end{eqnarray*}
\item  ii) If $a=(d+2q+2)^{-1}$, $\tilde a=(d+2q)^{-1}$, and if there exist 
$c,\tilde c\geq 0$ such that $\lim_{n \to \infty}{nh_n^{d+2q+2}}=c$ and 
$\lim_{n \to \infty}{n\tilde h_n^{d+2q}}=\tilde c$, 
then
\begin{eqnarray*}
\left(
\begin{array}{c}
\sqrt{nh_n^{d+2}}(\theta_n-\theta)\\
\sqrt{n\tilde{h}_n^{d}}(\tilde{\mu}_n-\mu)
\end{array}\right)\stackrel{\mathcal{D}}{\longrightarrow}                             \mathcal{N} 
\left(D(c,\tilde{c})AB_q(\theta), A\Sigma A\right).
\end{eqnarray*}
\item iii) If (C2) is satisfied, then
\begin{eqnarray*}
\left(
\begin{array}{c}
\frac{1}{h_n^q}(\theta_n-\theta)\\
\frac{1}{\tilde{h}_n^q}(\tilde{\mu}_n-\mu)
\end{array}\right)\stackrel{\mathbb{P}}{\longrightarrow}AB_q(\theta).
\end{eqnarray*}
\end{description}
\end{theo}

\begin{Rk}
\label{rque Rosenblatt}
The simultaneous weak convergence rate of nonrecursive estimators of the location and 
size of the mode can be established by following the lines of the proof 
of Theorem \ref{cen}. More precisely, set 
$$B_q^*(\theta) =  \left(
\begin{array}{c}
\frac{(-1)^{q}}{q!}\nabla\Big(\sum_{j=1}^d\beta_j^q\frac{\partial^qf}{\partial{x_j^q}}(\theta)\Big)\\
\frac{(-1)^q}{q!}\sum_{j=1}^d\beta_j^q\frac{\partial^qf}{\partial{x_j^q}}(\theta)
\end{array}\right),
\ 
\Sigma^* = \left(
\begin{array}{cc}
{f(\theta)G} & 0 \\
0 & {f(\theta)\int_{\mathbb{R}^d}K^2(z)dz}
\end{array}\right),
$$
let $\theta_n^*$ be Parzen's kernel estimator 
of the location of the mode and $\tilde\mu^*_n=\tilde f_n^*(\theta^*_n)$ be the 
kernel estimator of the size of the mode defined with the help of $\theta_n^*$ and 
of Rosenblatt's density estimator $\tilde f_n^*$ (the bandwidth $(\tilde h_n)$ 
defining $\tilde f_n^*$ being eventually different from the banwidth $(h_n)$ 
used to define $\theta_n^*$); Theorem \ref{cen} holds when 
$\theta_n$,  $\tilde\mu_n$, $B_q(\theta)$, $\Sigma$ are replaced by 
$\theta_n^*$,  $\tilde\mu^*_n$, $B_q^*(\theta)$, $\Sigma^*$, respectively. 
\end{Rk}

Part 1 and Part 2 in the case $c=\tilde c=0$ (respectively Part 3) of Theorem \ref{cen} correspond to the case when the bias (respectively the variances) of both estimators $\theta_n$ and $\tilde\mu_n$ are negligeable in front of their respective variances (respectively bias). When $c,\tilde c>0$, Part 2 of Theorem \ref{cen} corresponds to the case when the bias and the variance of each estimator $\theta_n$ and $\tilde \mu_n$ have the same convergence rate. Other possible conditions lead to different combinations; these ones have been omitted for sake of simplicity.\\

Theorem \ref{cen} gives the joint weak convergence rate of $\theta_n$ and 
$\tilde \mu_n$. Of course, it is also possible to estimate the location and 
the size of the mode separately. Concerning the estimation of the location of 
the mode, let us enlighten that the advantage of the semirecursive estimator 
$\theta_n$ on its nonrecursive version $\theta_n^*$ is that its asymptotic 
variance $[1+a(d+2)]^{-1}f(\theta)G$ is smaller than the one of Parzen's 
estimator, which equals $f(\theta)G$ (see, e.g. Romano \cite{romanom} for 
the case $d=1$ and Mokkadem and Pelletier \cite{mokkadempelletierm} for 
the case $d\geq 1$); this advantage of semirecursive estimators will be 
discussed again in Section \ref{2.3}. 
The estimation of the size of the mode is of course not independent of the 
estimation of the location, since the estimator $\tilde \mu_n$ is constructed 
with the help of the estimator $\theta_n$. To get a good estimation of 
the size of the mode, it seems obvious that $\theta_n$ should be computed 
with a bandwidth $(h_n)$ leading to its optimal convergence rate (or, at 
least, to a convergence rate close to the optimal one). The main information 
given by Theorem \ref{cen} is that, for $\tilde \mu_n$ to converge at the 
optimal rate, the use of a second bandwidth $(\tilde h_n)$ is then necessary.\\

Let us enlighten that, in the case when $\theta_n$ and $\tilde \mu_n$ 
satisfy a central limit theorem (Parts 1 and 2 of Theorem \ref{cen}), these 
estimators are asymptotically independent, although, in its definition, 
the estimator of the size of the mode is heavily connected to the one 
of the location of the mode. As pointed out by a referee, this property 
was expected. As a matter of fact (and as mentioned in the introduction), 
the location of the mode is a parameter which gives information on the shape 
of the density derivative, whereas the size of the mode gives information 
on the shape of the density itself. This constatation must be related to 
the fact that the weak (and strong) convergence rate of $\theta_n$ 
is given by the one of the gradient of $f_n$, whereas the 
 weak (and strong) convergence rate of $\tilde \mu_n$ is given by the one 
of $\tilde f_n$ itself; the variance of the density estimators converging to 
zero faster than the one of the estimators of the density derivatives, 
the asymptotic independence of $\theta_n$ and $\tilde \mu_n$ is completely 
explained. \\

Let us finally say one word on our assumptions on the bandwidths. In the 
framework of nonrecursive estimation, there is no need to assume that 
$(h_n)$ and $(\tilde h_n)$ are regularly varying sequences. In the case of 
semirecursive estimation, this assumption can obviously not be omitted, since 
the exponents $a$ and $\tilde a$ stand in the expressions of the asymptotic 
bias $B_q(\theta)$ and variance $\Sigma$. This might be seen as a slight 
inconvenient of semirecursive estimation; however, as it is enlightened in 
the following section, it turns out to be an advantage, since the asymptotic 
variances of the semirecursive estimators are smaller than the ones of 
the nonrecursive estimators.

\subsection{Construction of confidence regions and simulations studies}
\label{2.3}

The application of Theorem \ref{cen} (and of Remark \ref{rque Rosenblatt}) 
allows the construction of confidence regions (simultaneous or not) 
of the location and of the size of the mode, as well as 
confidence ellipsoids of the couple $(\theta,\mu)$.
Hall \cite{hall3m} shows that, in order to construct confidence regions, 
avoiding bias estimation by a slight undersmoothing is more efficient than 
explicit bias correction. In the framework of undersmoothing, the asymptotic bias 
of the estimator is negligeable in front of its asymptotic variance; 
according to the estimation by confidence regions point of view, the parameter 
to minimize is thus the asymptotic variance. Now, note that 
\begin{eqnarray*}
\Sigma = 
 \left(
\begin{array}{cc}
\left[1+a(d+2)\right]^{-1}I_d & 0 \\
0 & \left[1+\tilde{a}d\right]^{-1}
\end{array}\right)\Sigma^*
\end{eqnarray*}
(where $A\Sigma A$ (respectively $A\Sigma^* A$) is the asymptotic covariance matrix of 
the semirecursive estimators $(\theta_n,\tilde \mu_n)$ 
(respectively of the nonrecursive estimators $(\theta_n^*,\tilde \mu_n^*)$).
In order to construct confidence regions for the location and/or size of the mode, 
it is thus much preferable to use semirecursive estimators rather than nonrecursive 
estimators. Simulations studies confirm this theoritical conclusion, whatever the 
parameter ($\theta$, $\mu$ or $(\theta,\mu)$) for which confidence regions 
are contructed is. For sake of succintness, we do not give all these simulations
results here, but focuse on the construction of confidence ellipsoid for 
$(\theta,\mu)$; 
the aim of this example is of course to enlighten the advantage of using semirecursive 
estimators rather than nonrecursive estimators, but also to show how this confidence 
region gives informations on the shape of the density, and, consequently allows 
to measure the pertinence of the parameter location of the mode.
\\

To construct confidence regions for $(\theta,\mu)$, we consider the case $d=1$. 
The following corollary is a straightforward consequence of Theorem \ref{cen}.

\begin{cor}
\label{corollaire}
Let $\theta_n$ and $\tilde\mu_n$ be defined by (\ref{babamode}) and (\ref{bis}), 
respectively, and assume that (A1)-(A5) hold. Moreover, let $(h_n)$ and 
$(\tilde h_n)$ either satisfy (C1) or be such that   
$\lim_{n \to \infty}{nh_n^{2q+3}}=0$ and 
$\lim_{n \to \infty}{n\tilde h_n^{2q+1}}=0$ with $a=(2q+3)^{-1}$
and $\tilde a=(2q+1)^{-1}$.
We then have
\begin{eqnarray}
\label{convergence cor}
\frac{(1+3a){nh_n^3}[f''(\theta)]^2}{{ {f}(\theta)\int_{\mathbb{R}} K'^2(x)dx}}
(\theta_n-\theta)^2+
\frac{(1+\tilde a){n\tilde h_n}}{{f}(\theta)\int_{\mathbb{R}}K^2(x)dx}
(\tilde{\mu}_n-\mu)^2 & \stackrel{\mathcal{D}}{\longrightarrow} & \chi^2(2).
\end{eqnarray}
Moreover, (\ref{convergence cor}) still holds when the parameters ${f}(\theta)$ and 
$f''(\theta)$ are replaced by consistent estimators. 
\end{cor}

\begin{Rk}
\label{rque corollaire}
In view of Remark \ref{rque Rosenblatt}, in the case when the nonrecursive estimators 
$\theta_n^*$ and $\tilde\mu^*_n$ are used, 
(\ref{convergence cor}) becomes 
\begin{eqnarray}
\label{convergence rque cor}
\frac{{nh_n^3}[f''(\theta)]^2}{{ {f}(\theta)\int_{\mathbb{R}} K'^2(x)dx}}
(\theta_n^*-\theta)^2+
\frac{{n\tilde h_n}}{{f}(\theta)\int_{\mathbb{R}}K^2(x)dx}
(\tilde{\mu}_n^*-\mu)^2 & \stackrel{\mathcal{D}}{\longrightarrow} & \chi^2(2)
\end{eqnarray}
(and, again, this convergence still holds when the parameters ${f}(\theta)$ and 
$f''(\theta)$ are replaced by consistent estimators). 
\end{Rk}

Let $\check{f}_n''$ (respectively $\check{f}_n^{*''}$) be the recursive estimator 
(respectively the nonrecursive Rosenblatt's estimator) of $f''$ 
computed with the help of a bandwidth $\check{h}_n$, and set 
\begin{eqnarray*}
P_n=\frac{(1+3a)nh_n^3[\check{f}_n''(\theta_n)]^2}{\tilde{f}_n(\theta_n)
\int_{\mathbb{R}}K'^2(x)dx}, & & 
Q_n = \frac{(1+\tilde a){n\tilde h_n}}{{\tilde{f}_n(\theta_n)
\int_{\mathbb{R}}K^2(x)dx}}, \\
P_n^*=\frac{nh_n^3[\check{f}_n^{*''}(\theta_n^*)]^2}{\tilde{f}_n^*(\theta_n^*)
\int_{\mathbb{R}}K'^2(x)dx}, 
& & 
Q_n^* = \frac{{n\tilde h_n}}{{\tilde{f}^*_n(\theta_n^*)
\int_{\mathbb{R}}K^2(x)dx}}.
\end{eqnarray*}
Moreover, let $c_\alpha$ be such that $\mathbb{P}(Z\leq c_\alpha)=1-\alpha$, where 
$Z$ is $\chi^2(2)$-distributed; in view of Corollary~\ref{corollaire} and 
Remark \ref{rque corollaire}, the sets 
\begin{eqnarray*}
\mathcal{E}_{\alpha}& = &\left\{(\theta,\mu)/\  
P_n(\theta_n-\theta)^2+Q_n(\tilde{\mu}_n-\mu)^2 \leq c_{\alpha} \right\}\\
\mathcal{E}_{\alpha}^*& = &\left\{(\theta,\mu)/\  
P_n^*(\theta_n^*-\theta)^2+Q_n^*(\tilde{\mu}_n^*-\mu)^2 \leq c_{\alpha} \right\}
\end{eqnarray*}
are confidence ellipsoids for $(\theta,\mu)$ with asymptotic coverage level 
$1-\alpha$. Let us dwell on the fact that both confidence regions have the same 
asymptotic level, but the lengths of the axes of the first one (constructed with 
the help of the semirecursive estimators $\theta_n$ and $\tilde{\mu}_n$) are smaller 
than the ones of the second one (constructed with 
the help of the nonrecursive estimators $\theta_n^*$ and $\tilde{\mu}_n^*$). \\

We now present simulations results. In order to see the relationship between the shape 
of the confidence ellipsoids and the one of the density, the density $f$ we consider 
is the density of the ${\cal N}(0,\sigma^2)$-distribution, the parameter $\sigma$ 
taking the values $0.3$, $0.4$, $0.5$, $0.7$, $0.75$, $1$, $1.5$, $2$, and $2.5$. 
We use the sample size $n=100$ and the coverage level $1-\alpha=95\%$ 
(and thus $c_{\alpha}=5.99$). In each case, the number of simulations is $N=5000$. 
The kernel we use is the standard Gaussian density; the bandwidths are 
$$h_n=\dfrac{n^{-1/7}}{(\log n)},\ \ \tilde h_n=\dfrac{n^{-1/5}}{(\log n)},
\ \ \check{h}_n=n^{-1/9}.$$
Table 1 below gives, for each value of $\sigma$, the empirical values of 
$\theta_n$ , $\theta_n^*$, $\mu_n$ , $\mu_n^*$ (with respect to the 5000 
simulations), and:
\begin{description}
\item[$b$] the empirical length of the $\theta$-axis of the confidence ellipsoid 
$\mathcal{E}_{5\%}$;
\item[$b^*$] the empirical length of the $\theta$-axis of the confidence ellipsoid 
$\mathcal{E}_{5\%}^*$;
\item[$a$] the empirical length of  the $\mu$-axis of the confidence ellipsoid 
$\mathcal{E}_{5\%}$;
\item[$a^*$] the empirical length of the $\mu$-axis of the confidence ellipsoid 
$\mathcal{E}_{5\%}^*$;
\item[$p$] the empirical coverage level of the confidence ellipsoid 
$\mathcal{E}_{5\%}$;
\item[$p^*$] the empirical coverage level of the confidence ellipsoid 
$\mathcal{E}_{5\%}^*$.
\end{description}
\begin{center}
{\bf Table 1}\\
\begin{tabular}{lccccccccc}
$\sigma $ & $0.3$ & $0.4$ & $0.5$ & $0.7$ & $0.75$ & $1$ & $1.5$ & $2$ & $2.5
$ \\ \hline 
$\theta_n$ & $-0.002$ & $0.004$ & $0.001$ & $0.003$ & $0.002$ & $0.014$ & $-0.005$ & 
$-0.009$ & $0.014$ \\ 
$\theta_n^*$ & $0.003$ & $0.005$ & $0.001$ & $0.005$ & $-0.008$ & $0.016$ & $0.003$ & 
$-0.020$ & $-0.046$ \\ 
$b$ & $1.154$ & $1.346$ & $1.805$ & $2.898$ & $3.160$ & $5.218$ & $10.094$
& $17.866$ & $17.405$ \\ 
$b^*$ & $1.166$ & $1.458$ & $1.968$ & $3.300$ & $3.582$ & $5.925$ & $12.943$
& $21.946$ & $23.715$ \\ \hline 
$\mu_n$ & $1.335$ & $0.989$ & $0.782$ & $0.564$ & $0.522$ & $0.401$ & $0.263$ & 
$0.196$ & $0.155$ \\ 
$\mu_n^*$ & $1.312$ & $0.979$ & $0.783$ & $0.562$ & $0.512$ & $0.388$ & $0.269$ & 
$0.193$ & $0.163$ \\ 
$a$ & $0.444$ & $0.399$ & $0.365$ & $0.322$ & $0.315$ & $0.283$ & $0.247$ & 
$0.224$ & $0.210$ \\ 
$a^*$ & $0.514$ & $0.459$ & $0.420$ & $0.369$ & $0.363$ & $0.327$ & $0.287$ & 
$0.261$ & $0.246$ \\ \hline 
$p$ & $98.7\%$ & $97.8\%$ & $98.2\%$ & $98.4\%$ & $97.7\%$ & $97.8\%$ & $%
97.5\%$ & $97.2\%$ & $98.4\%$\\
$p^*$ & $98.6\%$ & $98.1\%$ & $98.4\%$ & $98.2\%$ & $96.8\%$ & $96.6\%$ & $%
96.9\%$ & $97.7\%$ & $98.2\%$  
\end{tabular}
\end{center}

Confirming our theoritical results, we see that the empirical coverage 
levels of both confidence ellipsoids 
$\mathcal{E}_{5\%}^*$ and $\mathcal{E}_{5\%}$ are similar, but that 
the empirical areas of the ellipsoids $\mathcal{E}_{5\%}$
(constructed with the help of the semirecursive estimators) are always 
smaller than the ones of the the ellipsoids $\mathcal{E}_{5\%}^*$ (constructed with 
the help of the nonrecursive estimators). \\

Let us now discuss the interest of the estimation of the size of the mode and 
the one of the joint estimation of the location and size of the mode. Both estimations 
give informations on the shape of the probability density and, consequently, 
allow to measure the pertinence of the parameter location of the mode. Of course, 
the parameter $\theta$ is significant only in the case when the high of the peak 
is large enough; since we consider here the example of the 
${\cal N}(0,\sigma^2)$-distribution, this corresponds to the case when $\sigma$ 
is small enough. Estimating only the size of the mode gives a first idea of the 
shape of the density around the location of the mode (for instance, when the
size is estimated around $0.16$, it is clear that the density is very flat). 
Now, the shape of the confidence ellipsoids allows to get a more precise idea. 
As a matter of fact, for small values of $\sigma$, the length of the $\mu$-axis
is larger than the one of the $\theta$-axis; as $\sigma$ increases, the length 
of the $\mu$-axis decreases, and the one of the $\theta$-axis increases 
(for $\sigma=2.5$, the length of the $\theta$-axis is larger than 20 times the one 
of the $\mu$-axis). Let us underline that these variations of the lengths of the 
axes are not due to bad estimations results; Table 2 below gives the values 
of the lengths $b$ (respectively $b^*$) of the $\theta$-axis, $a$ (respectively $a^*$) 
of the $\mu$-axis of the ellipsoids computed with the semirecursive estimators 
$\theta_n$ and $\tilde{\mu}_n$ (respectively with 
the nonrecursive estimators $\theta_n^*$ and $\tilde{\mu}_n^*$) in the case when 
the true values of the parameters ${f}(\theta)$ and $f''(\theta)$ are used 
(that is, by straightforwardly applying (\ref{convergence cor}) and 
(\ref{convergence rque cor})).
\\


\begin{center}
{\bf Table 2}\\
\begin{tabular}{lccccccccc}
$\sigma $ & $0.3$ & $0.4$ & $0.5$ & $0.7$ & $0.75$ & $1$ & $1.5$ & $2$ & $2.5
$ \\ \hline 
$b$ & $0.159$ & $0.327$ & $0.571$ & $1.357$ & $1.572$ & $3.227$ & $8.895$ & 
$18.260$ & $31.899$ \\ 
$b^*$ & $0.190$ & $0.390$ & $0.682$ & $1.622$ & $1.879$ & $3.858$ & $10.631$
& $21.825$ & $38.127$ \\ \hline 
$\mu$ & $1.333$ & $0.998$ & $0.798$ & $0.570$ & $0.532$ & $0.399$ & $0.266$ & 
$0.199$ & $0.159$ \\ 
$a$ & $0.465$ & $0.403$ & $0.360$ & $0.303$ & $0.294$ & $0.255$ & $0.208$ & 
$0.180$ & $0.161$ \\
$a^*$ & $0.509$ & $0.441$ & $0.395$ & $0.332$ & $0.322$ & $0.279$ & $0.228$ & 
$0.197$ & $0.176$  
\end{tabular}
\end{center}






\subsection{Strong convergence rate}
\label{2.4}

To establish the joint strong convergence rate of $\theta_n$ and 
$\tilde{\mu}_n$, we need the following additionnal assumption.
\begin{description}
\item (A6) i) $h$ and $\tilde{h}$ are differentiables, their derivatives vary regularly with exponent $(-a-1)$ and $(-\tilde{a}-1)$ respectively.\\
         ii) There exists $n_0\in \mathbb{N}$ such that
\begin{eqnarray*} 
n\geq m\geq n_0 \quad \Rightarrow \quad \max\left\{\frac{mh_m^{-(d+2)}}{nh_n^{-(d+2)}} ; \frac{m\tilde{h}_m^{-d}}{n\tilde{h}_n^{-d}}\right\} & = & 
 \frac{\min\left\{mh_m^{-(d+2)};m\tilde{h}_m^{-d}\right\}}{\min\left\{nh_n^{-(d+2)};n\tilde{h}_n^{-d}\right\}}. 
\end{eqnarray*}
\end{description}
\begin{Rk}\label{unex}
Assumption (A6)ii) holds when $a\neq \tilde a$, and in the case $a=\tilde a$, it is satisfied when $\mathcal{L}_{\theta}(n)=(\mathcal{L}_{\mu}(n))^{\frac{d}{d+2}}$ for $n$ large enough.
\end{Rk}
Moreover, condition (C2) is replaced by the following one.
\begin{description}
\item (C'2) Either (C2) i) is fulfilled or $\tilde a = \dfrac{1}{d+2q}$, 
$\lim_{n \to \infty}\dfrac{(\mathcal{L}_{\mu}(n))^{d+2q}}{2\log\log n}=\infty$, and 
$\dfrac{1}{2(d+2q)}<a <  \dfrac{1}{d+2q+2}$.
\end{description}

Before stating the almost sure convergence rate of $(\theta_n^T,\tilde{\mu}_n)^T$, let us remark that Proposition 2.3 in Mokkadem and Pelletier \cite{mokkadempelletierm} ensures that the matrix $G$ (and thus the matrix $\Sigma$) is nonsingular.
\begin{theo}\label{ite}
Let $\theta_n$ and $\tilde\mu_n$ be defined by (\ref{babamode}) and (\ref{bis}), respectively,
and assume that (A1)-(A6) hold. 
\begin{description}
\item i) If (C1) 
is fulfilled, then, with probability one, the sequence 
\begin{eqnarray*}
\frac{1}{\sqrt{2\log \log n}}\left(
\begin{array}{c}
\sqrt{nh_n^{d+2}}(\theta_n-\theta) \\
\sqrt{n\tilde{h}_n^{d}}(\tilde{\mu}_n-\mu)
\end{array}\right)
\end{eqnarray*}
is relatively compact and its limit set is the ellipsoid
\begin{eqnarray*}
\mathcal{E}=\left\{\nu \in \mathbb{R}^{d+1}\ \mbox{such that} \ \nu^TA^{-1}\Sigma^{-1}A^{-1}\nu\leq 1 \right\}.
\end{eqnarray*}
\item ii) If $a=(d+2q+2)^{-1}$, $\tilde a=(d+2q)^{-1}$, and if there exist 
$c, \tilde c\geq 0$ such that \\ $\lim_{n \to \infty}{nh_n^{d+2q+2}}
/{(2\log\log n)}=c$ and  
$\lim_{n \to \infty}{n\tilde h_n^{d+2q}}/{(2\log\log n)}=\tilde c$, 
then, 
with probability one, the sequence 
\begin{eqnarray*}
\frac{1}{\sqrt{2\log \log n}}\left(
\begin{array}{c}
\sqrt{nh_n^{d+2}}(\theta_n-\theta) \\
\sqrt{n\tilde{h}_n^{d}}(\tilde{\mu}_n-\mu)
\end{array}\right)
\end{eqnarray*}
is relatively compact and its limit set is the ellipsoid
$$\displaystyle\mathcal{E}=\left\{\nu \in \mathbb{R}^{d+1}\ \mbox{such that} \ \left(A^{-1}\nu -D\left(c,\tilde{c}\right)B_q(\theta)\right)^T\Sigma^{-1}\left(A^{-1}\nu-D\left(c,\tilde{c}\right)B_q(\theta)\right)\leq 1 \right\}.$$
\item iii) If (C'2) is satisfied, then
\begin{eqnarray*}
\left(
\begin{array}{c}
\frac{1}{h_n^q}(\theta_n-\theta)\\
\frac{1}{\tilde{h}_n^q}(\tilde{\mu}_n-\mu)
\end{array}\right)\stackrel{\mbox{a.s.}}{\longrightarrow}AB_q(\theta).
\end{eqnarray*}
\end{description}
\end{theo}
\begin{Rk}
(C'1) implies that 
$\lim_{n \to \infty}{nh_n^{d+2q+2}}/{\log\log n}=0$ and 
$\lim_{n \to \infty}{n\tilde h_n^{d+2q}}/{\log\log n}=0$, whereas 
(C'2) implies that 
$\lim_{n \to \infty}{nh_n^{d+2q+2}}/{\log\log n}=\infty$ and 
$\lim_{n \to \infty}{n\tilde h_n^{d+2q}}/{\log\log n}=\infty$.
\end{Rk}

Laws of the iterated logarithm for Parzen's nonrecursive kernel mode 
estimator were established by Mokkadem and Pelletier \cite{mokkadempelletierm}. 
The technics of demonstration used in the framework of nonrecursive 
estimators are totally different from those employed to prove 
Theorem \ref{ite}. This is due to the following fondamental difference 
between the nonrecursive estimator $\theta_n^*$ and the semirecursive 
estimator $\theta_n$: the study of the asymptotic behaviour of $\theta_n^*$ 
comes down to the one of a \emph{triangular sum} of independent variables, 
whereas the study of the asymptotic behaviour of $\theta_n$ reduces to the 
one of a \emph{sum} of independent variables. Of course, this difference is 
not quite important for the study of the weak convergence rate. But, for 
the study of the strong convergence rate, it makes the case of the 
semirecursive estimation much easier than the case of the nonrecursive 
estimation. In particular, on the oppposite to the weak convergence rate, 
the joint strong convergence rate of the nonrecursive estimators 
$\theta_n^*$ and $\tilde\mu_n^*$ cannot be obtained by following the lines 
of the proof of Theorem \ref{ite}, and remains an open question.

\section{Proofs}

Let us first note that an important consequence of (A3) which will be used throughout the proofs 
is that 
\begin{eqnarray}
\label{relmode}
\mbox{if} \ \ \beta a <1,\ \ \mbox{then} \ \ \lim_{n \to \infty}\frac{1}{nh_n^{\beta}}\sum_{i=1}^nh_i^{\beta} & = &  \frac{1}{1-a\beta}.
\end{eqnarray}
Moreover, for all $\varepsilon>0$ small enough,
\begin{eqnarray}
\label{cremode}
\frac{1}{n}\sum_{i=1}^nh_i^q=O\left(h_n^{q-\varepsilon}+\frac{1}{n}\right).
\end{eqnarray}
As a matter of fact: (i) 
if $a q<1$, \eqref{cremode} follows easily from \eqref{relmode}; (ii)
if $a q>1$, since $\sum_ih_i^q$ is summable, \eqref{cremode} holds; (iii)
if $a q=1$, since $a (q-\varepsilon)<1$, using \eqref{relmode} again, we have 
$n^{-1}\sum_{i=1}^nh_i^q=O(h_n^{q-\varepsilon})$, and thus \eqref{cremode} follows.
Of course \eqref{relmode} and \eqref{cremode} also hold when $(h_n)$ and $a$ are replaced by 
$(\tilde h_n)$ and $\tilde a$, respectively.\\

Our proofs are now organized as follows. Section \ref{4.1} is devoted to the proof of the strong 
consistency of $\theta_n$ and $\tilde{\mu}_n$.
In Section \ref{4.2}, we give the convergence rate of the derivatives of $f_n$. 
In Section \ref{4.3}, we show how the study of the joint weak and strong convergence rate 
of $\theta_n$ and  $\tilde{\mu}_n$ can be related to the one of $\nabla f_n(\theta)$ and 
$\tilde{f}_n(\theta)$. In Section \ref{4.4} (respectively in Section \ref{4.5}), 
we establish the joint weak convergence rate (respectively the joint strong convergence rate) 
of $\nabla f_n(\theta)$ and $\tilde{f}_n(\theta)$.  
Finally, Section \ref{4.6} is devoted to the proof of Theorems \ref{cen} and \ref{ite}.

\subsection{Proof of Proposition \ref{Prop}} \label{4.1}
 Since $\theta_n$ is the mode of $f_n$ and $\theta$ the mode of $f$, we have:
\begin{eqnarray}
\label{bobomode}
0  \leq f(\theta)-f(\theta_n) & = &[f(\theta)-f_n(\theta_n)]+[f_n(\theta_n)-f(\theta_n)]
\leq [f(\theta)-f_n(\theta)]+[f_n(\theta_n)-f(\theta_n)] \nonumber \\
& \leq & \big|f(\theta)-f_n(\theta)\big|+ \big|f_n(\theta_n)-f(\theta_n)\big| 
\leq 2\|f_n-f\|_{\infty}.  
\end{eqnarray}
The application of Theorem 5 in Mokkadem et al. \cite{mokkadempelletierthiamm} with $|\alpha| =0$ and $v_n=\log n$ ensures that for any $\delta>0$, there exists $c(\delta)>0$ such that 
$\mathbb{P}[(\log n) \|f_n-\mathbb{E}(f_n)\|_{\infty} \geq \delta] \leq 
\exp(-c(\delta){\sum_{i=1}^nh_i^d}/{(\log n)^2})$.
In view of \eqref{relmode}, since $ad<1$, we can write
$$
n^2\exp\left(-c(\delta)\frac{\sum_{i=1}^nh_i^d}{(\log n)^2}\right)  =  n^2\exp\left(-c(\delta)\frac{nh_n^d}{(\log n)^2}\frac{\sum_{i=1}^nh_i^d}{nh_n^d}\right)  =  o(1).
$$
Borell-Cantelli's Lemma ensures that 
$\lim_{n \to \infty}\|f_n-\mathbb{E}(f_n)\|_{\infty}=0$ a.s. 
Since $\lim_{n \to \infty}\|\mathbb{E}(f_n)-f\|_{\infty}=0$, it follows from 
\eqref{bobomode} that 
$\lim_{n \to \infty}f(\theta_n)=f(\theta)$ a.s.
Since $f$ is continuous, $\lim_{\|z\| \to \infty}f(z)=0$ and $\theta$ is the unique mode of $f$, 
we deduce that $\lim_{n \to \infty}\theta_n=\theta$ a.s. 
Now, we have 
$$|\tilde{\mu}_n-\mu| \leq 
 |\tilde{f}_n(\theta_n)-f(\theta_n)|+|f(\theta_n)-f(\theta)|  \leq 
 \|\tilde{f}_n-f\|_{\infty} +2\|f_n-f\|_{\infty},$$
where the last inequality follows from (\ref{bobomode}). As previously, one can show that $\lim_{n \to \infty} \|\tilde f_n-f\|_{\infty}=0$ and thus $\lim_{n \to \infty}\tilde{\mu}_n=\mu$ a.s.

\subsection{Convergence rate of the derivatives of the density}\label{4.2}
 For any $d$-uplet $\displaystyle[\alpha]=\Big(\alpha_1,\ldots,\alpha_d\Big)\in \mathbb{N}^d$, we set $|\alpha|=\alpha_1+\cdots+\alpha_d$ and, for any function $g$, let 
$\partial^{[\alpha]}g(x)={\partial^{|\alpha|}g}/({\partial x_1^{\alpha_1}\ldots\partial x_d^{\alpha_d}})(x)$
denote the $[\alpha]$-th partial derivative of $g$. 

\begin{lemma}\label{aea} 
Assume (A3)-(A5) hold. Let $(g_n)$ and $(b_n)$ be defined as follows:
\begin{eqnarray}
\label{sonemode}
\left\{
\begin{array}{ll}
g_n=f_n \ \ \mbox{and} \ \ b_n=h_n \ \ \mbox{or}\\
g_n=\tilde {f}_n \ \ \mbox{and} \ \ b_n=\tilde {h}_n.
\end{array}
\right.
\end{eqnarray}
For $|\alpha|\in\{0,1,2\}$, we have
\begin{eqnarray*}
\lim_{n \to \infty}\frac{n}{\sum_{i=1}^nb_i^q}\bigg[\mathbb{E}\big[\partial^{[\alpha]}g_n(x)\big]-\partial^{[\alpha]}f(x)\bigg]  & = & \frac{(-1)^q}{q!}\partial^{[\alpha]}\left(\sum_{j=1}^d\beta_j^q\frac{\partial^qf}{\partial x_j^q}\right)(x)
\end{eqnarray*}
where $\beta_j^q$ is defined in (\ref{recmode}).
Moreover, if we set $M_q = \sup_{x \in \mathbb{R}^d}\|D^q\partial^{[\alpha]}f(x)\|$, then
\begin{eqnarray*}
\lim_{n \to \infty}\frac{n}{\sum_{i=1}^nb_i^q}\sup_{x \in \mathbb{R}^d}\left|\mathbb{E}\left(\partial^{[\alpha]}g_n(x)\right)-\partial^{[\alpha]}f(x)\right|  & \leq & \frac{M_q}{q!}\int_{\mathbb{R}^d}\|z\|^q\left|K(z)\right|dz.
\end{eqnarray*} 
\end{lemma}

\begin{lemma}\label{aia}
Let $U$ be a compact set of $\mathbb{R}^d$ and assume that (A1)iii), (A3), (A4) and (A5)ii) hold. Let $(g_n)$  and $(b_n)$ be defined as in \eqref{sonemode}. Then, 
for all $\gamma >0$ and $|\alpha|=1,2$, we have
\begin{eqnarray*}
\sup_{x \in U}\left|\partial^{[\alpha]}g_n(x) -\mathbb{E}\left(\partial^{[\alpha]}g_n(x)\right)\right| & = &  O\left(\sqrt{\frac{(\log n)^{1+\gamma}}{\sum_{i=1}^nb_i^{d+2|\alpha|}}}\right) \ \ \mbox{ a.s.}\\
\end{eqnarray*}
\end{lemma}
Lemma \ref{aea} is proved in Mokkadem et al. \cite{mokkadempelletierthiamm}. We now prove 
Lemma \ref{aia}. Set 
$v_n  =  [\sum_{i=1}^nb_i^{d+2|\alpha|}]^{1/2}$ $[(\log n)^{1+\gamma}]^{-1/2}$.
Applying  Proposition 3 in Mokkadem et al. \cite{mokkadempelletierthiamm}, it holds that for any $\delta>0$, there exists $c(\delta)>0$ such that 
\begin{eqnarray*}
\mathbb{P}\left[\sup_{x \in U}v_n\left|\partial^{[\alpha]}g_n(x)-\mathbb{E}\left(\partial^{[\alpha]}g_n(x)\right)\right| \geq \delta\right] & \leq & 
\exp\left(-c(\delta)\frac{\sum_{i=1}^nb_i^{d+2|\alpha|}}{2v_n^2}\right).
\end{eqnarray*}
Since $\lim_{n\to \infty}\sum_{i=1}^nb_i^{d+2|\alpha|}/(v_n^2\log n)=\infty $ we have, for $n$ large enough, 
$c(\delta){\sum_{i=1}^nb_i^{d+2|\alpha|}}/({2v_n^2})\geq  2\log n$, 
and Lemma \ref{aia} follows from the application of Borel-Cantelli's Lemma.
\subsection{Relationship between  
\boldmath$((\theta_n-\theta)^T,(\tilde{\mu}_n-\mu))^T$ 
and \boldmath${( [\nabla f_n(\theta)]^T,\tilde{f}_n(\theta)-f(\theta))}^T$} \label{4.3}
By definition of $\theta_n$, we have $\nabla f_n(\theta_n)=0$ so that 
\begin{eqnarray}
\label{totomode}
\nabla f_n(\theta_n)-\nabla f_n(\theta)= -\nabla f_n(\theta).
\end{eqnarray}
For each $i \in \left\{1, \ldots,d\right\}$, a Taylor expansion applied to the real valued application 
${\partial f_n}/{\partial x_i}$ implies the existence of $\varepsilon_n(i)=(\varepsilon_n^{(1)}(i), \dots, \varepsilon_n^{(d)}(i))^t$ such that 
\begin{eqnarray*}
\left \lbrace
\begin{array}{l}
\frac{\partial f_n}{\partial x_i}(\theta_n)-\frac{\partial f_n}{\partial x_i}(\theta)=\sum_{j=1}^d\frac{\partial^2f_n}{\partial x_i\partial x_j}(\varepsilon_n(i))\Big(\theta_n^{(j)}-\theta^{(j)}\Big),\\
\big|\varepsilon_n^{(j)}(i)-\theta^{(j)}\big|\leq\big|\theta_n^{(j)}(i)-\theta^{(j)}\big| \quad  \forall j \in \left\{1, \dots, d\right\}.\\
\end{array}
\right.
\end{eqnarray*}
Define the $d\times d$ matrix $H_n=(H_n^{(i,j)})_{1\leq i,j\leq d}$ by setting 
$H_n^{(i,j)}=\frac{\partial^2 f_n}{\partial x_i\partial x_j}(\varepsilon_n(i))$;
Equation (\ref{totomode}) can be then  rewritten as 
$H_n(\theta_n-\theta)=-\nabla f_n(\theta)$.
Now, set 
\begin{eqnarray}
R_n & = & \tilde{f}_n(\theta_n)-\tilde{f}_n(\theta).
\label{tetemode}
\end{eqnarray}  
We can then write:
\begin{eqnarray}
\label{tatamode}
\left(
\begin{array}{c}
\left[D^2f(\theta)\right]^{-1}H_n(\theta_n-\theta)\\
\tilde{\mu}_n-\mu
\end{array}\right) & = & 
\left(
\begin{array}{c}
-\left[D^2f(\theta)\right]^{-1}\nabla f_n(\theta)\\
\tilde{f}_n(\theta)-f(\theta)
\end{array}\right)   +  
\left(                                                                                                         \begin{array}{c}
0 \\
R_n
\end{array}\right).
\end{eqnarray}
Let $U$ be a compact set of $\mathbb{R}^d$ containing $\theta$. The combination of Lemmas \ref{aea} 
and \ref{aia} with $|\alpha|=2$, $g_n=f_n$ and $b_n=h_n$ ensures that for any $\gamma>0$ and $\varepsilon >0$ small enough,
\begin{eqnarray}
\label{diammode}
\sup_{x \in U}\left|\partial^{[\alpha]}f_n(x) -\partial^{[\alpha]}f(x)\right| & = & O\left(\sqrt{\frac{(\log n)^{1+\gamma}}{\sum_{i=1}^nh_i^{d+4}}}+\frac{\sum_{i=1}^nh_i^q}{n}\right) \ \  \mbox{a.s.} \nonumber\\
& = & O\left(\sqrt{\frac{(\log n)^{1+\gamma}}{nh_n^{d+4}}}+h_n^{q-\varepsilon}+\frac{1}{n}\right)
=o(1)\ \ \mbox{a.s.} 
\end{eqnarray}
Since $D^2f$ is continuous in a neighbourhood of $\theta$ and since $\lim_{n  \to \infty}\theta_n=\theta$ a.s., \eqref{diammode} ensures that $\lim_{n \to \infty}H_n=D^2f(\theta)$ a.s. It follows that the weak and a.s. behaviours of $((\theta_n-\theta)^T,(\tilde{\mu}_n-\mu))^T$ are given by the one of the right-hand-sided term of (\ref{tatamode}).

\subsection{Weak convergence rate of \boldmath${( [\nabla f_n(\theta)]^T,
\tilde{f}_n(\theta)-f(\theta))}^T$} \label{4.4}
Let us at first assume that the following lemma holds.

\begin{lemma}\label{ada}
Let Assumptions (A1)i), (A1)iv), (A3), (A4)i) and (A4)ii) hold. Then
 \begin{eqnarray*}
W_n & = &\left(
\begin{array}{c}
\sqrt{nh_n^{d+2}}\Big[\nabla f_n(\theta)-\mathbb{E}\big(\nabla f_n(\theta)\big)\Big]\\
\sqrt{n\tilde{h}_n^{d}}\Big[\tilde{f}_n(\theta)-\mathbb{E}\big(\tilde{f}_n(\theta)\big)\Big]
\end{array}\right)\stackrel{\mathcal{D}}{\longrightarrow}\mathcal{N}\Big(0,\Sigma\Big).
\end{eqnarray*}
\end{lemma}
The application of Lemma \ref{aea} gives 
\begin{eqnarray}
\label{fatmode}
\lim_{n \to \infty}
\left(
\begin{array}{c}
\frac{n}{\sum_{i=1}^nh_i^q}\mathbb{E}\big(\nabla f_n(\theta)\big)\\
\frac{n}{\sum_{i=1}^n\tilde{h}_i^q}\Big[\mathbb{E}\big(\tilde{f}_n(\theta)\big)-f(\theta)\Big] 
\end{array}
\right) & = & 
 \left(
\begin{array}{c}
\frac{(-1)^{q}}{q!}\nabla\Big(\sum_{j=1}^d\beta_j^q\frac{\partial^qf}{\partial{x_j^q}}(\theta)\Big)\\
\frac{(-1)^q}{q!}\sum_{j=1}^d\beta_j^q\frac{\partial^qf}{\partial{x_j^q}}(\theta)
\end{array}\right). 
\end{eqnarray}

\begin{description}
\item 1) If $aq< 1$ and $\tilde aq< 1$, by  using \eqref{relmode}, it is straightforward to see that
\begin{eqnarray}
\label{fatemode}
\lim_{n \to \infty}
\left(
\begin{array}{c}
\frac{1}{h_n^q}\mathbb{E}\big(\nabla f_n(\theta)\big)\\
\frac{1}{\tilde{h}_n^q}\Big[\mathbb{E}\big(\tilde{f}_n(\theta)\big)-f(\theta)\Big] 
\end{array}
\right)  & = & B_q(\theta). 
\end{eqnarray}
\item 2) Let us now consider the case $aq\geq1$ and $\tilde aq\geq1$. We have
\begin{eqnarray*}
\sqrt{nh_n^{d+2}}\mathbb{E}\left(\nabla f_n(\theta)\right) & = &\sqrt{nh_n^{d+2}} \frac{\sum_{i=1}^nh_i^q}{n}\frac{n}{\sum_{i=1}^nh_i^q}\mathbb{E}\left(\nabla f_n(\theta)\right),
\end{eqnarray*}
with, in view of \eqref{cremode}, for all $\varepsilon >0$ small enough,
$$\sqrt{nh_n^{d+2}}\frac{\sum_{i=1}^nh_i^q}{n}  =  O\left(n^{\frac{1}{2}(1-(a-\varepsilon)(d+2))}n^{-aq+a\varepsilon}\right) 
=o(1).$$
Applying \eqref{fatmode}, it follows that 
$\lim_{n \to \infty}\sqrt{nh_n^{d+2}}\mathbb{E}(\nabla f_n(\theta))= 0$.
Proceeding in the same way for $\mathbb{E}(\tilde f_n(\theta))$, we obtain 
\begin{eqnarray}
\label{dermode}
\lim_{n \to \infty}\left(
\begin{array}{c}
\sqrt{nh_n^{d+2}}\mathbb{E}\big(\nabla f_n(\theta)\big)\\
\sqrt{n\tilde{h}_n^{d}}\Big[\mathbb{E}\big(\tilde{f}_n(\theta)\big)-f(\theta)\Big]
\end{array}\right)& = & 0.
\end{eqnarray}
\end{description}
The combination of either \eqref{fatemode} or \eqref{dermode} and of Lemma \ref{ada} gives the weak convergence rate of \\${( [\nabla f_n(\theta)]^T,\tilde{f}_n(\theta)-f(\theta))}^T$:
\begin{itemize}
\item If (C1) holds, then 
\begin{eqnarray}
\label{cramode}
\left(
\begin{array}{c}
\sqrt{nh_n^{d+2}}\nabla f_n(\theta)\\
\sqrt{n\tilde{h}_n^{d}}(\tilde{f}_n(\theta)-f(\theta))
\end{array}\right)\stackrel{\mathcal{D}}{\longrightarrow}                             \mathcal{N}\left( 0,
 \Sigma\right).
\end{eqnarray}
\item If $a=(d+2q+2)^{-1}$, $\tilde a=(d+2q)^{-1}$, and if there exist 
$c,\tilde c\geq 0$ such that $\lim_{n \to \infty}{nh_n^{d+2q+2}}=c$ and 
$\lim_{n \to \infty}{n\tilde h_n^{d+2q}}=\tilde c$, then 
\begin{eqnarray}
\label{crabismode}
\left(
\begin{array}{c}
\sqrt{nh_n^{d+2}}\nabla f_n(\theta)\\
\sqrt{n\tilde{h}_n^{d}}(\tilde{f}_n(\theta)-f(\theta))
\end{array}\right)\stackrel{\mathcal{D}}{\longrightarrow}                             \mathcal{N}\left( D(c,\tilde{c})B_q(\theta),
 \Sigma\right).
\end{eqnarray}
\item  If (C2) holds, since $aq<1$ and $\tilde aq<1$, \eqref{relmode} implies that
\begin{eqnarray}
\label{cromode}
\left(
\begin{array}{c}
\frac{1}{h_n^q}\nabla f_n(\theta)\\
\frac{1}{\tilde{h}_n^q}(\tilde{f}_n(\theta)-f(\theta))
\end{array}\right)\stackrel{\mathbb{P}}{\longrightarrow} B_q(\theta).
\end{eqnarray}
\end{itemize} 
\paragraph{Proof of Lemma \ref{ada}}
To prove Lemma \ref{ada}, we first prove that
\begin{eqnarray}
\label{papamode}
\lim_{n \to \infty}\mathbb{E}\Big(W_nW_n^T\Big) =  \Sigma, 
\end{eqnarray}
and then check that $(W_n)$ satisfies Lyapounov's condition. Set
\begin{eqnarray*}
Y_{k,n} & = & \frac{1}{\sqrt{nh_n^{-d-2}}}h_k^{-d-1}\left[\nabla K\left(\frac{\theta-X_k}{h_k}\right)-\mathbb{E}\left(\nabla K\left(\frac{\theta-X_k}{h_k}\right)\right)\right]\\
Z_{k,n} & = & \frac{1}{\sqrt{n\tilde{h}_n^{-d}}}\tilde{h}_k^{-d}\left[K\left(\frac{\theta-X_k}{\tilde{h}_k}\right)-\mathbb{E}\left(K\left(\frac{\theta-X_k}{\tilde{h}_k}\right)\right)\right],
\end{eqnarray*}
and note that 
\begin{eqnarray*}
\mathbb{E}\Big(W_nW_n^T\Big) & = & \sum_{k=1}^n\left(
\begin{array}{cc}
\mathbb{E}\Big(Y_{k,n}Y_{k,n}^T\Big) & \mathbb{E}\Big(Y_{k,n}Z_{k,n}\Big) \\
\mathbb{E}\Big(Y_{k,n}^TZ_{k,n}\Big) & \mathbb{E}\Big(Z_{k,n}^2\Big) 
\end{array}\right).
\end{eqnarray*}
Now, for any $\displaystyle s,t \in \left\{1,\ldots,d\right\}$, we have
\begin{eqnarray*}
 \mathbb{E}\Bigg[\frac{\partial K}{\partial x_s}\left(\frac{\theta-X_k}{h_k}\right)\frac{\partial K}{\partial x_t}\left(\frac{\theta-X_k}{h_k}\right)\Bigg] 
& = & \int_{\mathbb{R}^d}\frac{\partial K}{\partial x_s}\left(\frac{\theta-y}{h_k}\right)\frac{\partial K}{\partial x_t}\left(\frac{\theta-y}{h_k}\right)f(y)dy  \\ 
& = & h_k^df(\theta)G_{s,t}+o(h_k^d), 
\end{eqnarray*}
and since,
$\mathbb{E}\left[\frac{\partial K}{\partial x_s}\left(\frac{\theta-X_k}{h_k}\right)\right] = O(h_k^d)$, we deduce that 
\begin{eqnarray}
\label{ricmode}
&\mathbb{E}\Bigg(\Bigg[\nabla K\left(\frac{\theta-X_k}{h_k}\right)-\mathbb{E}\Bigg(\nabla K\left(\frac{\theta-X_k}{h_k}\right)\Bigg)\Bigg]\Bigg[\nabla K\left(\frac{\theta-X_k}{h_k}\right)-\mathbb{E}\Bigg(\nabla K\left(\frac{\theta-X_k}{h_k}\right)\Bigg)\Bigg]^T\Bigg)  &\nonumber \\ &  =   f(\theta)Gh_k^d\Big[1+o(1)\Big] & 
\end{eqnarray}
which implies that 
$\lim_{n \to \infty}\sum_{k=1}^n\mathbb{E}(Y_{k,n}Y_{k,n}^T)={f(\theta)}[1+a(d+2)]^{-1}G$.
In the same way, we have  
\begin{eqnarray}
\label{rocmode}
\mathbb{E}\Bigg(\Bigg[K\left(\frac{\theta-X_k}{\tilde{h}_k}\right)-\mathbb{E}\Bigg(K\left(\frac{\theta-X_k}{\tilde{h}_k}\right)\Bigg)\Bigg]^2\Bigg) & = &  \tilde{h}_k^df(\theta)\int_{\mathbb{R}^d}K^2(z)dz\Big[1+o(1)\Big]
\end{eqnarray}
and thus $\lim_{n \to \infty}\sum_{k=1}^n\mathbb{E}(Z_{k,n}^2)={f(\theta)}[1+\tilde ad]^{-1}
\int_{\mathbb{R}^d}K^2(z)dz$. 
Moreover, set $h_n^*=\min(h_n,\tilde{h}_n)$; we have
\begin{eqnarray*}
\mathbb{E}\Bigg[\nabla K\left(\frac{\theta-X_k}{h_k}\right)K\left(\frac{\theta-X_k}{\tilde{h}_k}\right)\Bigg] & = &
h_k^{*d}\int_{\mathbb{R}^d}\nabla K\Big(\frac{h_k^*}{h_k}z\Big)K\Big(\frac{h_k^*}{\tilde{h}_k}z\Big)f(\theta-h_k^*z)dz. 
\end{eqnarray*}
Noting that $f(\theta-h_k^*z)=f(\theta)+h_k^*R_k(\theta,z)$ with $\left|R_k(\theta,z)\right| \leq \|\nabla f\|_{\infty}\|z\|$, we get
\begin{eqnarray*}
\lefteqn{\mathbb{E}\Bigg[\nabla K\left(\frac{\theta-X_k}{h_k}\right)K\left(\frac{\theta-X_k}{\tilde{h}_k}\right)\Bigg]}\\ & = & h_k^{*d}\left[f(\theta)\int_{\mathbb{R}^d}\nabla K\Big(\frac{h_k^*}{h_k}z\Big)K\Big(\frac{h_k^*}{\tilde{h}_k}z\Big)dz+h_k^*\int_{\mathbb{R}^d}\nabla K\Big(\frac{h_k^*}{h_k}z\Big)K\Big(\frac{h_k^*}{\tilde{h}_k}z\Big)R_k(\theta,z)dz \right].
\end{eqnarray*}
Since the function $z\mapsto \big[\nabla K(z)\big]K(z$ is odd (in each coordinate), the first right-handed integral is zero, and, since $h_k^*$ equals either $h_k$ or $\tilde h_k$, we get                 
\begin{eqnarray*}
\lefteqn{\left|\left|\mathbb{E}\Bigg[\nabla K\left(\frac{\theta-X_k}{h_k}\right)K\left(\frac{\theta-X_k}{\tilde{h}_k}\right)\Bigg]\right|\right|}\\  & \leq & h_k^{*(d+1)}\|\nabla f\|_{\infty}\left[\|K\|_{\infty}\int_{\mathbb{R}^d}\|z\|\|\nabla K(z)\|dz+\|\nabla K\|_{\infty}\int_{\mathbb{R}^d}\|z\||K(z)|dz\right] 
=O\left(h_k^{*(d+1)}\right).
\end{eqnarray*} 
We then deduce that
\begin{eqnarray}
\label{rcemode}
&&\mathbb{E}\Bigg(\Bigg[\nabla K\left(\frac{\theta-X_k}{h_k}\right)-\mathbb{E}\Bigg(\nabla K\left(\frac{\theta-X_k}{h_k}\right)\Bigg)\Bigg]\Bigg[ K\left(\frac{\theta-X_k}{\tilde{h}_k}\right)-\mathbb{E}\Bigg(K\left(\frac{\theta-X_k}{\tilde{h}_k}\right)\Bigg)\Bigg]\Bigg) \nonumber\\   & = &
  O\left(\left[\min(h_k,\tilde{h}_k)\right]^{d+1}\right)+O\left(h_k^d\tilde h_k^d\right)
= 
 O\left(h_k^{\frac{d+1}{2}}\tilde{h}_k^{\frac{d+1}{2}}\right),
\end{eqnarray}
and thus, in view of \eqref{relmode},
$$\sum_{k=1}^n\mathbb{E}\Big(Y_{k,n}Z_{k,n}\Big) =  O\left(\frac{1}{\sqrt{(nh_n^{-d-2})(n\tilde{h}_n^{-d})}}\sum_{k=1}^nh_k^{-\frac{d+1}{2}}\tilde{h}_k^{\frac{1-d}{2}}\right)
= o(1),$$
which concludes the proof of (\ref{papamode}).
Now we check that $(W_n)$ satisfies the Lyapounov's condition. Set $p>2$. Since $K$ and $\nabla K$ are bounded and integrable, we have $\int_{\mathbb{R}^d}\|\nabla K(z)\|^pdz<\infty$ and $\int_{\mathbb{R}^d}|K(z)|^pdz<\infty$. 
It follows that 
\begin{eqnarray*}
\sum_{k=1}^n\mathbb{E}\left(\|Y_{k,n}\|^p\right) & = &  
O\left(\frac{1}{(nh_n^{-d-2})^{\frac{p}{2}}}\sum_{k=1}^nh_k^{(-d-1)p}\int_{\mathbb{R}^d}\left|\left|\nabla K\left(\frac{\theta-y}{h_k}\right)\right|\right|^pf(y)dy\right) \nonumber \\
& = & 
O\Bigg(\frac{1}{(nh_n^{-d-2})^{\frac{p}{2}}}\sum_{k=1}^nh_k^{(-d-1)p}h_k^d\Bigg) 
 = 
 o(1), \\
\sum_{k=1}^n\mathbb{E}\left(\big|Z_{k,n}\big|^p\right) & = & O\left(\frac{1}{(n\tilde{h}_n^{-d})^{\frac{p}{2}}}\sum_{k=1}^n\tilde{h}_k^{-dp}\int_{\mathbb{R}^d}\Big|K\left(\frac{\theta-y}{\tilde{h}_k}\right)\Big|^pf(y)dy\right) \nonumber \\ & = & O\Bigg(\frac{1}{(n\tilde{h}_n^{-d})^{\frac{p}{2}}}\sum_{k=1}^n\tilde{h}_k^{-dp}\tilde{h}_k^d\Bigg) 
=o(1),
\end{eqnarray*} 
which concludes the proof of Lemma \ref{ada}.  

\subsection{A.s.  convergence rate of \boldmath${( [\nabla f_n(\theta)]^T,\tilde{f}_n(\theta)-f(\theta))}^T$} \label{4.5}
Let us at first assume that the following lemma holds.
\begin{lemma}\label{map}
Let Assumptions (A1)i), (A1)iv), (A3), (A4)i), (A4)ii) and (A6) hold. With probability one, the sequence
 \begin{eqnarray*}
 \frac{1}{\sqrt{2\log \log n}}\left(
\begin{array}{c}
\sqrt{nh_n^{d+2}}\left[\nabla f_n(\theta)-\mathbb{E}\big(\nabla f_n(\theta)\big)\right]\\
\sqrt{n\tilde{h}_n^{d}}\left[\tilde{f}_n(\theta)-\mathbb{E}\big(\tilde{f}_n(\theta)\big)\right]
\end{array}\right)
\end{eqnarray*}
is relatively compact and its limit set is $\mathcal{E}=\left\{\nu \in \mathbb{R}^{d+1} \ \mbox{such  that} \ \nu^T\Sigma^{-1}\nu\leq 1 \right\}$.
\end{lemma}
The combination of either \eqref{fatemode} or \eqref{dermode} and of Lemma \ref{map} gives the almost sure convergence rate of \\${( [\nabla f_n(\theta)]^T,\tilde{f}_n(\theta)-f(\theta))}^T$:  
\begin{itemize}
\item If (C1) holds, then, with probability one, the sequence 
\begin{eqnarray}
\label{oulmode}
\frac{1}{\sqrt{2\log \log n}}\left(
\begin{array}{c}
\sqrt{nh_n^{d+2}}\nabla f_n(\theta)\\
\sqrt{n\tilde{h}_n^{d}}\left[\tilde{f}_n(\theta)-f(\theta)\right]
\end{array}\right)
\end{eqnarray}
is relatively compact and its limit set is 
$\mathcal{E}=\left\{\nu \in \mathbb{R}^{d+1} \ \mbox{such that} \ \nu^T\Sigma^{-1} \nu \leq 1 \right\}$.
\item If $a=(d+2q+2)^{-1}$, $\tilde a=(d+2q)^{-1}$, and if there exist 
$c,\tilde c\geq 0$ such that $\lim_{n \to \infty}{nh_n^{d+2q+2}}/$ ${(2\log\log n)}=c$ and  $\lim_{n \to \infty}{n\tilde h_n^{d+2q}}/{(2\log\log n)}=\tilde c$,  
then with probability one, the sequence 
\begin{eqnarray}
\label{oulbismode}
\frac{1}{\sqrt{2\log \log n}}\left(
\begin{array}{c}
\sqrt{nh_n^{d+2}}\nabla f_n(\theta)\\
\sqrt{n\tilde{h}_n^{d}}\left[\tilde{f}_n(\theta)-f(\theta)\right]
\end{array}\right)
\end{eqnarray}
is relatively compact and its limit set is 
\begin{eqnarray*}
\mathcal{E}=\left\{\nu \in \mathbb{R}^{d+1} \ \mbox{such that} \ \left(\nu-D(c,\tilde{c})B_q(\theta)\right)^T\Sigma^{-1} \left(\nu-D(c,\tilde{c})B_q(\theta)\right) \leq 1 \right\}.
\end{eqnarray*} 
\item If (C'2) holds, then 
\begin{eqnarray}
\label{madomode}
\left(
\begin{array}{c}
\frac{1}{h_n^q}\nabla f_n(\theta)\\
\frac{1}{\tilde{h}_n^q}[\tilde{f}_n(\theta)-f(\theta)]
\end{array}\right)\stackrel{\mbox{a.s.}}{\longrightarrow} B_q(\theta).
\end{eqnarray}
\end{itemize}
We now prove Lemma \ref{map}.
Set
$$\Gamma = f(\theta)\left(
\begin{array}{cc}
G & 0 \\
0 & \int_{\mathbb{R}^d}K^2(z)dz
\end{array}\right),\ \ 
\Delta_n  =  \left(
\begin{array}{cc}
\frac{1}{\sqrt{nh_n^{-d-2}}}I_d & 0\\
0 & \frac{1}{\sqrt{n\tilde{h}_n^{-d}}}
\end{array}\right),\ \  Q_n = \left(
\begin{array}{cc}
\sqrt{h_n^{-d-2}}I_d & 0\\
0 & \sqrt{\tilde{h}_n^{-d}}
\end{array}\right),$$ 
let $(\varepsilon_n)$ be a sequence of $\mathbb{R}^{d+1}$-valued, independent and $\mathcal{N}(0, \Gamma)$-distributed random vectors, and set 
$S_n =\sum_{k=1}^nQ_k\varepsilon_k$. 
In order to prove Lemma \ref{map}, we first establish the following Lemma 
\ref{aha} in Section \ref{Sec a}, and then show in Section \ref{Sec b} how 
Lemma \ref{map} can be deduced from Lemma \ref{aha}.
\begin{lemma}\label{aha}
Let Assumptions (A1)i), (A1)iv), (A3), (A4)i), (A4)ii) and (A6)ii) hold.
With probability one, the sequence $(T_n)\equiv({\Sigma^{-1/2}\Delta_nS_n}/
{\sqrt{2\log \log n}})$ is relatively compact and its limit set is the unit ball  $\displaystyle\overline{\mathcal{B}}_{d+1}(0,1)=\left\{\nu \in \mathbb{R}^{d+1} \ \mbox{such  that} \  \|\nu\|_2 \leq 1 \right\}$.
\end{lemma}

\subsubsection{Proof of Lemma \ref{aha}} 
\label{Sec a}
Set $B_n = \mathbb{E}(S_nS_n^T)$,
let $\|x\|_2$ (respectively ${\left|\!\left|\!\left|A\right|\!\right|\!\right|_2}$) denote the euclidean norm (respectively the spectral norm) of the vector $x$ (respectively of the matrix $A$).
The application of Theorem 2 in Koval \cite{kovalm} ensures that
\begin{eqnarray*}
\limsup_{n \to \infty}\frac{\|\Sigma^{-1/2}\Delta_nS_n\|_2}{\sqrt{2\left|\!\left|\!\left|\Sigma^{-1/2}\Delta_nB_n\Delta_n\Sigma^{-1/2}\right|\!\right|\!\right|_2\log\log \left|\!\left|\!\left|B_n\right|\!\right|\!\right|_2}} \leq 1 \quad \mbox{a.s.}
\end{eqnarray*}
Since $\lim_{n \to \infty}\Delta_nB_n\Delta_n=\Sigma$ and $\log\log \left|\!\left|\!\left|B_n\right|\!\right|\!\right|_2 \sim \log\log n$, we deduce that 
\begin{eqnarray}
\label{zizimode}
\limsup_{n \to \infty}\|T_n\|_2\leq 1 \quad \mbox{a.s.}
\end{eqnarray}
Thus, the sequence $(T_n)$ is relatively compact and its limit set $\mathcal{U}$ is included in $\overline{\mathcal{B}}_{d+1}(0,1)$.
Now, set $\displaystyle \mathcal{S}_{d+1}=\left\{w \in \mathbb{R}^{d+1},\ \|w\|_2=1 \right\}$, and let us at first assume that 
\begin{eqnarray}
\label{tapmode}
\forall w \in  \mathcal{S}_{d+1}, \ \ \limsup_{n \to \infty}w^TT_n\geq 1 \quad \mbox{a.s.}
\end{eqnarray}
The combination of (\ref{zizimode}) and (\ref{tapmode}) ensures that, with probability one,
$\displaystyle\forall \varepsilon >0, \forall n_0\geq 1,\exists n\geq n_0$ such  that $\displaystyle w^TT_n>1-\varepsilon$ and $\|T_n\|_2^2\leq 1+\varepsilon$.  
Noting that $\displaystyle \|T_n-w\|_2^2=\|T_n\|_2^2+\|w\|_2^2-2w^TT_n$, it follows that, with probability one, $\forall \varepsilon >0$, $\forall n_0\geq 1$, 
$\exists n\geq n_0$ such that 
$\|T_n-w\|_2^2\leq 1+\varepsilon+1-2(1-\varepsilon)=3\varepsilon$.
Thus, with probability one, $\mathcal{S}_{d+1}\subset\mathcal{U}$. To deduce that $\overline{\mathcal{B}}_{d+1}(0,1) \subset \mathcal{U}$, we introduce $(e_k)$, a sequence of real-valued, independent, and $\mathcal{N}(0,1)$-distributed random variables such that $(e_k)$ is independent of $(\varepsilon_k)$. Moreover, we set 
\begin{eqnarray*}
\tilde{Q}_n=\left(
\begin{array}{cc}
\sqrt{h_n^{-d-2}}I_{d+1} & 0\\
0 & \sqrt{\tilde{h}_n^{-d}}
\end{array}\right),  \quad  \tilde{S}_n & = & \sum_{k=1}^n\tilde{Q}_k\left(
\begin{array}{c}
e_k\\
\varepsilon_k
\end{array}\right)
\end{eqnarray*}
\begin{eqnarray*}
\tilde{\Delta}_n=
 \left(
\begin{array}{cc}
\frac{1}{\sqrt{nh_n^{-d-2}}}I_{d+1} & 0\\
0 & \frac{1}{\sqrt{n\tilde{h}_n^{-d}}}
\end{array}\right)
, \ \ \mbox{and} \ \ \tilde{\Sigma}=\left(
\begin{array}{cc}
1 & 0 \\
0 & \Sigma
\end{array}\right).
\end{eqnarray*}
We then note that the previous result applied to $(\tilde{T}_n)\equiv
({\tilde\Sigma^{-1/2}\tilde{\Delta}_n\tilde{S}_n}/{\sqrt{2\log \log n}})$  
ensures that, with probability one, $\mathcal{S}_{d+2}=\{w \in \mathbb{R}^{d+2},\ \|w\|_2=1\}$ is included in the limit set of $\tilde{T}_n$. Now let $\pi : \mathbb{R}^{d+2}\longrightarrow \mathbb{R}^{d+1}$ be the projection map defined by $\pi((x_1, \ldots, x_{d+2})^T)=(x_2, \ldots, x_{d+2})^T$. We clearly have $\pi(\mathcal{S}_{d+2})=\overline{\mathcal{B}}_{d+1}(0,1)$ 
and $\pi(\tilde{T}_n)=T_n$, and thus deduce that, with probability one, $\overline{\mathcal{B}}_{d+1}(0,1)$ is included in the limit set of $T_n$. 
To conclude the proof of Lemma \ref{aha}, it remains to prove (\ref{tapmode}). In fact, we shall prove that, 
\begin{eqnarray}
\label{eqamode}
\forall w\neq 0, \ \ \limsup_{n \to \infty}\frac{w^T\Delta_nS_n}{\sqrt{2\log \log n}} \geq \sqrt{w^T\Sigma w}\quad \mbox{a.s.}
\end{eqnarray}
Set $v_n= \min\{[nh_n^{-(d+2)}]^{1/2}; [n\tilde{h}_n^{-d}]^{1/2}\}$,
$\displaystyle A_n =v_nw^T\Delta_n$ and $\displaystyle V_n=\mathbb{E}\big(A_nS_nS_n^TA_n^T\big)$; we follow a method used by Petrov \cite{petrovm} in the proof of his Theorems 7.1 and 7.2. 
Since $\displaystyle\lim_{n \to \infty} V_n = \infty$, $\forall \tau >0$, there exists a non-decreasing sequence of integers ${n_k}$ such that $n_k\rightarrow\infty$ as $k\rightarrow \infty$ and $V_{n_k-1}\leq(1+\tau)^k\leq V_{n_k}$, $(k=1,2,\ldots)$.
Since $\lim_{n \to \infty}{V_{n-1}}/{V_n}=1$, we obtain $V_{n_k}\sim (1+\tau)^k$.
Moreover, we have 
\begin{eqnarray}
\label{yeyemode}
V_{n_k}-V_{n_{k-1}} = V_{n_k}\big(1-\frac{V_{n_{k-1}}}{V_{n_k}}\big)\sim V_{n_k}\frac{\tau}{\tau +1}. 
\end{eqnarray}
Set
$$\chi(n)  =  \sqrt{2V_n\log \log V_n},\ \ 
\psi(n_k)  =  \sqrt{2(V_{n_k}-V_{n_{k-1}})\log \log(V_{n_k}-V_{n_{k-1}})}.$$
It follows from (\ref{yeyemode}) that $\psi(n_k)\sim \tau^{1/2}\chi(n_{k-1})$. Then for any $\gamma \in ]0,1[$ and $k$ sufficiently large, we have 
\begin{eqnarray}
\label{yuyumode}
\lefteqn{\mathbb{P}\left(A_{n_k}S_{n_k}-A_{n_k}S_{n_{k-1}} \geq (1-\gamma)\psi(n_k)\right)} \nonumber  \\  
& \geq &                                                         \mathbb{P}\left(A_{n_k}S_{n_k}\geq (1 - \frac{\gamma}{2})\psi(n_k)\right)-\mathbb{P}\left(A_{n_k}S_{n_{k-1}} \geq \frac{\gamma \psi(n_k)}{2}\right) \nonumber \\  & \geq &                                                             \mathbb{P}\left(A_{n_k}S_{n_k}\geq (1 - \frac{\gamma}{2})\chi(n_k)\right)-\mathbb{P}\left(A_{n_k}S_{n_{k-1}} \geq \frac{\gamma \sqrt{\tau}}{3}\chi(n_{k-1})\right).   
\end{eqnarray}
Since $A_{n_k}S_{n_k}$ is $\mathcal{N}\big(0,V_{n_{k}}\big)$-distributed, we have 
\begin{eqnarray}
\label{wowomode}
\mathbb{P}\left(A_{n_k}S_{n_k} \geq (1- \frac{\gamma}{2})\chi(n_k) \right) & = & \nonumber \frac{1}{\sqrt{2\pi}}\int_{(1-\frac{\gamma}{2})\sqrt{2V_{n_{k}}\log \log V_{n_{k}}}}^{\infty}\exp{\left(-\frac{t^2}{2}\right)}dt \\ &\geq & \left[\log V_{n_k}\right]^{-(1+\mu)(1-\frac{\gamma}{2})^2}
\end{eqnarray}
for every $\mu$ and sufficiently large $k$.
Set $\displaystyle\tilde{V}_{n_k}=v_{n_k}^2w^T\Delta_{n_k}B_{n_{k-1}}\Delta_{n_k}w$; since $A_{n_k}S_{n_{k-1}}$ is  $\mathcal{N}\big(0,\tilde{V}_{n_k}\big)$-distributed, we have 
\begin{eqnarray*}
\mathbb{P}\Big(A_{n_k}S_{n_{k-1}}  \geq  \frac{\gamma \sqrt{\tau}}{3}\chi(n_{k-1})\Big)  = 
 \frac{1}{\sqrt{2\pi}}\int_{\frac{\gamma \sqrt{\tau}}{3}\sqrt{2\frac{V_{n_{k-1}}}{\tilde{V}_{n_k}}\log \log V_{n_{k-1}}}}^{\infty}\exp{\left(-\frac{t^2}{2}\right)}dt. 
\end{eqnarray*}
Let $\rho_{min}(A)$ (respectively $\rho_{max}(A)$) denote the smallest (respectively the largest) eigenvalue of a matrix $A$, set $\Sigma_n=\Delta_nB_n\Delta_n$, and note that 
\begin{eqnarray}
\label{thesmode}
\frac{V_{n_{k-1}}}{\tilde{V}_{n_k}} 
& \geq & \frac{v_{n_{k-1}}^2\rho_{min}(\Sigma_{n_{k-1}})}{v_{n_{k}}^2\rho_{max}(\Delta_{n_{k}}\Delta_{n_{k-1}}^{-1}\Sigma_{n_{k-1}}\Delta_{n_{k-1}}^{-1}\Delta_{n_{k}})}
\end{eqnarray}
with 
\begin{eqnarray}
\label{thesemode}
\rho_{max}(\Delta_{n_{k}}\Delta_{n_{k-1}}^{-1}\Sigma_{n_{k-1}}\Delta_{n_{k-1}}^{-1}\Delta_{n_{k}}) & \leq & \left|\!\left|\!\left|\Sigma_{n_{k-1}}\Delta_{n_{k-1}}^{-1}\Delta_{n_{k}}\Delta_{n_{k}}\Delta_{n_{k-1}}^{-1}\right|\!\right|\!\right|_2 \nonumber \\ & \leq & \left|\!\left|\!\left|\Sigma_{n_{k-1}}\right|\!\right|\!\right|_2 \left|\!\left|\!\left|\Delta_{n_{k-1}}^{-1}\Delta_{n_{k}}\Delta_{n_{k}}\Delta_{n_{k-1}}^{-1}\right|\!\right|\!\right|_2.
\end{eqnarray}
It follows from \eqref{relmode} and Assumption A6)ii) that 
\begin{equation}
\label{thesesmode}
\left|\!\left|\!\left|\Delta_{n_{k-1}}^{-1}\Delta_{n_{k}}\Delta_{n_{k}}\Delta_{n_{k-1}}^{-1}\right|\!\right|\!\right|_2  =   \max\left\{\frac{n_{k-1}h_{n_{k-1}}^{-(d+2)}}{n_kh_{n_k}^{-(d+2)}} ; \frac{n_{k-1}\tilde{h}_{n_{k-1}}^{-d}}{n_k\tilde{h}_{n_{k}}^{-d}}\right\}  \sim  
\frac{v_{n_{k-1}}^2}{v_{n_k}^2}.
\end{equation}
>From \eqref{thesmode}, \eqref{thesemode} and \eqref{thesesmode}, we deduce that, for sufficiently large $k$,
$$
\frac{V_{n_{k-1}}}{\tilde{V}_{n_k}}  \geq  \frac{\rho_{min}(\Sigma_{n_{k-1}})}{2\rho_{max}(\Sigma_{n_{k-1}})}  \geq  \frac{\rho_{min}(\Sigma)}{4\rho_{max}(\Sigma)}
$$
and therefore, for sufficiently large $k$,
\begin{eqnarray}
\label{wiwimode}
\mathbb{P}\left(A_{n_k}S_{n_{k-1}} \geq \frac{\gamma \sqrt{\tau}}{3}\chi(n_{k-1})\right) \nonumber & \leq & 
\frac{1}{\sqrt{2\pi}}\int_{\frac{\gamma \sqrt{\tau}}{6}\sqrt{\frac{\rho_{min}(\Sigma)}{\rho_{max}(\Sigma)}}\sqrt{2V_{n_{k-1}}\log \log V_{n_{k-1}}}}^{\infty} \exp{\left(-\frac{t^2}{2}\right)}dt  \nonumber\\ & \leq  &
\left[\log V_{n_{k-1}}\right]^{\frac{-\gamma^2 \tau \rho_{min}(\Sigma)}{36 \rho_{max}(\Sigma)}}.
\end{eqnarray}
The inequalities (\ref{yuyumode}), (\ref{wowomode}) and (\ref{wiwimode}) imply that 
\begin{eqnarray*}
\mathbb{P}\left(A_{n_k}S_{n_k}-A_{n_k}S_{n_{k-1}} \geq (1-\gamma)\psi(n_k)\right)  \geq
\left[\log V_{n_k}\right]^{-(1+\mu)(1-\frac{\gamma}{2})^2} -\left[\log V_{n_{k-1}}\right]^{\frac{-\gamma^2 \tau \rho_{min}(\Sigma)}{36 \rho_{max}(\Sigma)}}.
\end{eqnarray*}
Thus, for sufficiently large $k$ and $\tau$, there exists $c >0$ such that $c$ does not depend on $k$ and
\begin{eqnarray*}
\mathbb{P}\left(A_{n_k}S_{n_k}-A_{n_k}S_{n_{k-1}} \geq (1-\gamma)\psi(n_k)\right)   \geq 
c\Big[k^{{-(1+\mu)(1-\frac{\gamma}{2})^2}}-k^{-1}\Big].
\end{eqnarray*} 
Choosing $\mu$ such that $\left(1+\mu\right)\left(1-{\gamma}/{2}\right)^2<1$, we get
\begin{eqnarray*}
\mathbb{P}\left(A_{n_k}S_{n_k}-A_{n_k}S_{n_{k-1}}  \geq (1-\gamma)\psi(n_k)\right)  \geq \frac{c}{2}k^{{-(1+\mu)(1-\frac{\gamma}{2})^2}}  
\end{eqnarray*} 
and thus 
$\sum_{k}\mathbb{P}(A_{n_k}S_{n_k}-A_{n_k}S_{n_{k-1}} \geq (1-\gamma)\psi(n_k))
= \infty$. 
Applying Borel-Cantelli's Lemma, we obtain
\begin{eqnarray}
\label{zazamode}
\mathbb{P}\left(A_{n_k}S_{n_k}-A_{n_k}S_{n_{k-1}} \geq (1-\gamma)\psi(n_k)\quad \mbox{i.o.}\right)=1.
\end{eqnarray}
Now, 
\begin{eqnarray*}
\limsup_{k \to \infty}\frac{\left|A_{n_k}S_{n_{k-1}}\right|}{\chi(n_{k-1})} 
& \leq &
\limsup_{k \to \infty}\frac{v_{n_{k}}\|w\|_2 \left|\!\left|\!\left|\Delta_{n_{k}}\Delta_{n_{k-1}}^{-1}\right|\!\right|\!\right|_2 \|\Delta_{n_{k-1}}S_{n_{k-1}}\|_2}{\sqrt{2v_{n_{k-1}}^2(w^T\Delta_{n_{k-1}}B_{n_{k-1}}\Delta_{n_{k-1}}w)\log \log V_{n_{k-1}}}} \\ & \leq &
\limsup_{k \to \infty}\frac{v_{n_{k}}\|w\|_2 \left|\!\left|\!\left|\Delta_{n_{k}}\Delta_{n_{k-1}}^{-1}\right|\!\right|\!\right|_2 \|\Delta_{n_{k-1}}S_{n_{k-1}}\|_2}{\sqrt{2v_{n_{k-1}}^2(w^T\Sigma w)\log \log V_{n_{k-1}}}}.
\end{eqnarray*}
Applying  Theorem 2 in Koval \cite{kovalm} again, and using the fact that $\lim_{n \to \infty} \Delta_nB_n\Delta_n=\Sigma$, we obtain 
$\limsup_{n \to \infty}{\|\Delta_nS_n\|_2}/{\sqrt{2\left|\!\left|\!\left|\Sigma \right|\!\right|\!\right|_2\log \log n}} \leq 1$ a.s.
Therefore, 
\begin{eqnarray*}
\limsup_{k \to \infty}\frac{\left|A_{n_k}S_{n_{k-1}}\right|}{\chi(n_{k-1})} & \leq & \limsup_{k \to \infty}\frac{v_{n_{k}}\|w\|_2 \left|\!\left|\!\left|\Delta_{n_{k}}\Delta_{n_{k-1}}^{-1}\right|\!\right|\!\right|_2\sqrt{\left|\!\left|\!\left|\Sigma\right|\!\right|\!\right|_2}}{v_{n_{k-1}}\sqrt{w^T\Sigma w}}\quad\mbox{a.s.} 
\end{eqnarray*}
Since $ |\!|\!|\Delta_{n_{k}}\Delta_{n_{k-1}}^{-1}|\!|\!|_2 =
[\rho_{max}(\Delta_{n_{k-1}}^{-1}\Delta_{n_{k}}\Delta_{n_{k}}
\Delta_{n_{k-1}}^{-1})]^{1/2} \leq {2v_{n_{k-1}}}/{v_{n_{k}}}$, for sufficiently large $k$, we obtain 
$\limsup_{k \to \infty}{\left|A_{n_k}S_{n_{k-1}}\right|}/{\chi(n_{k-1})} \leq 
{2\|w\|_2\sqrt{\left|\!\left|\!\left|\Sigma\right|\!\right|\!\right|_2}}/{\sqrt{w^T\Sigma w}}$ a.s. 
Set $\varepsilon \in ]0,1[$ and \\ $\kappa={2\|w\|_2
\sqrt{\left|\!\left|\!\left|\Sigma\right|\!\right|\!\right|_2}}/{\sqrt{w^T\Sigma w}}$. 
Noting that
\begin{eqnarray*}
(1-\gamma)\psi(n_k)-2\kappa\chi(n_{k-1}) \sim 
\Big[(1-\gamma)\sqrt{\tau}(1+\tau)^{-1/2}-2\kappa(1+\tau)^{-1/2}\Big]\chi(n_k),
\end{eqnarray*}
and noting that $\gamma$ can be chosen sufficiently small and  $\tau$ sufficiently large so that $(1-\gamma)\sqrt{\tau}(1+\tau)^{-1/2}-2\kappa(1+\tau)^{-1/2} > 1-\varepsilon$,
we obtain
\begin{eqnarray*}
\mathbb{P}\left(A_{n_k}S_{n_k}>(1-\varepsilon)\chi(n_k)\quad \mbox{i.o.}\right)  \geq  
 \mathbb{P}\left(A_{n_k}S_{n_k} >(1-\gamma)\psi(n_k)-2\kappa\chi(n_{k-1})\quad \mbox{i.o.}\right).
\end{eqnarray*}
Taking (\ref{zazamode}) into account, we then obtain
$\mathbb{P}\left(A_{n_k}S_{n_k}>(1-\varepsilon)\chi(n_k)\mbox{ i.o.} \right) =1$.
We thus get\\ $\limsup_{n \to \infty} {A_nS_n}/{\chi(n)} \geq 1$ a.s.,
which proves (\ref{eqamode}), and concludes the proof of Lemma \ref{aha}.  

\subsubsection{Proof of Lemma \ref{map}} \label{Sec b}
Now, set
\begin{eqnarray*}
\tilde{V}_k & = & \left( \begin{array}{c}
h_k^{-d/2}\left[\nabla K\left(\frac{\theta-X_k}{h_k}\right)-\mathbb{E}\left(\nabla K\left(\frac{\theta-X_k}{h_k}\right)\right)\right]\\
\tilde{h}_k^{-d/2}\left[K\left(\frac{\theta-X_k}{\tilde{h}_k}\right)-\mathbb{E}\left(K\left(\frac{\theta-X_k}{\tilde{h}_k}\right)\right)\right]
\end{array}\right)
\end{eqnarray*}
and $\Gamma_k   = \mathbb{E}(\tilde{V}_k\tilde{V}_k^T)$. In view of (\ref{ricmode}), (\ref{rocmode}) and (\ref{rcemode}), we have $\lim_{k \to \infty}\Gamma_k=\Gamma$. It follows that $\exists k_0 \geq  1$ such that $\forall k\geq k_0$, $\Gamma_k$ is inversible; without loss of generality, we assume $k_0=1$, and set $\tilde{U}_k =\Gamma_k^{-1/2}\tilde{V}_k$.
Set $p \in ]2,4[$ and let $\mathcal{L}$ be a slowly varying function; we have:
\begin{eqnarray*}
\frac{\mathbb{E}\big(\|\tilde{U}_k\|^p\big)}{(k\log \log k)^{p/2}} & = &
O\left(\frac{h_k^{-dp/2}\mathbb{E}\left[\|\nabla K\left(\frac{\theta-X_k}{h_k}\right)\|^p\right]+
\tilde{h}_k^{-dp/2}\mathbb{E}\left[\left|K\left(\frac{\theta-X_k}{h_k}\right)\right|^p\right]
}{(k\log \log k)^{p/2}}\right)
\\ & =&
O\left(\frac{h_k^{d-dp/2}+
\tilde{h}_k^{d-dp/2}
}{(k\log \log k)^{p/2}}\right)
 \\ & = & 
O\left(\mathcal{L}(k)\left[k^{-\left[1+(\frac{p}{2}-1)(1-ad)\right]}+k^{-\left[1+(\frac{p}{2}-1)(1-\tilde{a}d)\right]}\right]\right)
\end{eqnarray*}
so that $\sum_k(k\log \log k)^{-p/2}\mathbb{E}(\|\tilde{U}_k\|^p) <\infty$.
By application of Theorem 2 of Einmahl \cite{einmahlm},
we deduce that  
$\sum_{k=1}^n\tilde{U}_k-\sum_{k=1}^n\eta_k=o(\sqrt{n\log \log n})$ a.s.,
where $\eta_k$ are independent, and $\mathcal{N}(0,I_{d+1})$-distributed random vectors.
It follows that 
\begin{eqnarray}
\label{wourimode}
\sum_{k=1}^n\Gamma^{1/2}\Gamma_k^{-1/2}\tilde{V}_k-\sum_{k=1}^n\varepsilon_k=o(\sqrt{n\log \log n})\quad \mbox{a.s.}
\end{eqnarray}
Now,
\begin{eqnarray*}
\lefteqn{\Delta_n\left[\sum_{k=1}^nQ_k\Gamma^{1/2}\Gamma_k^{-1/2}\tilde{V}_k - \sum_{k=1}^nQ_k\varepsilon_k\right]} \\ & = &
\Delta_n\sum_{k=1}^nQ_k\left[\Gamma^{1/2}\Gamma_k^{-1/2}\tilde{V}_k -\varepsilon_k\right] \\ & = & 
 \Delta_n\sum_{k=1}^nQ_k\left(\sum_{j=1}^k\left[\Gamma^{1/2}\Gamma_j^{-1/2}\tilde{V}_j -\varepsilon_j\right]-\sum_{j=1}^{k-1}\left[\Gamma^{1/2}\Gamma_j^{-1/2}\tilde{V}_j -\varepsilon_j\right]\right) \quad (\mbox{with} \sum_{j=1}^0=0)\\
& = &
\Delta_n\sum_{k=1}^{n-1}\left(Q_k-Q_{k+1}\right)\left(\sum_{j=1}^k\left(\Gamma^{1/2}\Gamma_j^{-1/2}\tilde{V}_j -\varepsilon_j\right)\right)+\Delta_nQ_n\sum_{j=1}^n\left(\Gamma^{1/2}\Gamma_j^{-1/2}\tilde{V}_j -\varepsilon_j\right) \\ & = & 
\Delta_n\sum_{k=1}^{n-1}\left(Q_k-Q_{k+1}\right)\left[o\left(\sqrt{k\log\log k}\right)\right]+\Delta_nQ_n\left[o\left(\sqrt{n\log\log n}\right)\right] \quad \mbox{a.s.}
\end{eqnarray*}
Moreover,
\begin{eqnarray*}
\lefteqn{\Delta_n\sum_{k=1}^{n-1}\left(Q_k-Q_{k+1}\right)\left[o\left(\sqrt{k\log\log k}\right)\right]} \\ & = &
\left(
\begin{array}{cc}
 \sqrt{\frac{h_n^{d+2}}{n}}\sum_{k=1}^{n-1}\left(h_k^{-\frac{d+2}{2}}-h_{k+1}^{-\frac{d+2}{2}}\right)o\left(\sqrt{k\log\log k}\right) & 0\\
0 &  \sqrt{\frac{\tilde{h}_n^d}{n}}\sum_{k=1}^{n-1}\left(\tilde{h}_k^{-\frac{d}{2}}-\tilde{h}_{k+1}^{-\frac{d}{2}}\right)
\end{array}\right)
 \\ & = & 
\left(
\begin{array}{cc}
o\left(\sqrt{h_n^{d+2}\log\log n}\right)\sum_{k=1}^{n-1}\left(h_k^{-\frac{d+2}{2}}-h_{k+1}^{-\frac{d+2}{2}}\right) & 0\\
0 & o\left(\sqrt{\tilde{h}_n^d\log\log n}\right)\sum_{k=1}^{n-1}\left(\tilde{h}_k^{-\frac{d}{2}}-\tilde{h}_{k+1}^{-\frac{d}{2}}\right)
\end{array}\right).
\end{eqnarray*}
Set $\phi(s)=\left[h(s)\right]^{-\frac{d+2}{2}}$ and $\tilde{\phi}(s)=\left[\tilde{h}(s)\right]^{-\frac{d}{2}}$, and let $u_k \in [k,k+1]$; since $\phi'$ and $\tilde{\phi}'$ vary  regularly with exponent 
$(a(d+2)/2-1)$ and $(\tilde{a}d/2-1)$ respectively, we have
$$\sum_{k=1}^{n-1}\left(h_k^{-\frac{d+2}{2}}-h_{k+1}^{-\frac{d+2}{2}}\right)  =  
O\left(\sum_{k=1}^{n-1}\phi'(u_k)\right) 
=O\left(\int_{1}^n\phi'(s)ds\right)  
 =  O\left(h_n^{-\frac{d+2}{2}}\right)$$
and 
$$
\sum_{k=1}^{n-1}\left(\tilde{h}_k^{-\frac{d}{2}}-\tilde{h}_{k+1}^{-\frac{d}{2}}\right)  =  
O\left(\sum_{k=1}^{n-1}\tilde{\phi}'(u_k)\right) 
=O\left(\int_{1}^n\tilde{\phi}'(s)ds\right) 
 =  O\left(\tilde{h}_n^{-\frac{d}{2}}\right),
$$
so that $\Delta_n\sum_{k=1}^{n-1}\left(Q_k-Q_{k+1}\right)\left[o\left(\sqrt{k\log\log k}\right)\right]  =  o\left(\sqrt{\log\log n}\right)$.
Since
$\Delta_nQ_n\left[o\left(\sqrt{n\log\log n}\right)\right] =o\left(\sqrt{\log\log n}\right)$,
we deduce that 
\begin{eqnarray*}
\frac{\Delta_n\sum_{k=1}^nQ_k\Gamma^{1/2}\Gamma_k^{-1/2}\tilde{V}_k}{\sqrt{2\log\log n}} - \frac{\Delta_n\sum_{k=1}^nQ_k\varepsilon_k}{\sqrt{2\log\log n}} = o(1) \quad \mbox{a.s.}
\end{eqnarray*}
The application of Lemma \ref{aha} then ensures that, with probability one, 
the sequence 
$(\Delta_n \sum_{k=1}^nQ_k\Gamma^{1/2}$ $\Gamma_k^{-1/2}\tilde{V}_k/
{\sqrt{2\log\log n}})$
is relatively compact and its limit set is $\mathcal{E}=\{\nu \in \mathbb{R}^{d+1} \ \mbox{such  that} \ \nu^T\Sigma^{-1}\nu \leq 1 \}$. Since
\begin{eqnarray*}
\frac{\Delta_n\sum_{k=1}^nQ_k\tilde{V}_k}{\sqrt{2\log\log n}} & = &
 \frac{\Delta_n\sum_{k=1}^nQ_k\Gamma^{1/2}\Gamma_k^{-1/2}\tilde{V}_k}{\sqrt{2\log\log n}}+
 \frac{\Delta_n\sum_{k=1}^nQ_k\left(I_{d+1}-\Gamma^{1/2}\Gamma_k^{-1/2}\right)\tilde{V}_k}{\sqrt{2\log\log n}} 
\end{eqnarray*}
with $\lim_{k \to \infty}(I_{d+1}-\Gamma^{1/2}\Gamma_k^{-1/2})=0$, 
Lemma \ref{map} follows. 
\subsection{
Proof of Theorems \ref{cen} and \ref{ite}} \label{4.6}
In view of (\ref{tatamode}) (and the comment below), Theorem \ref{cen} (respectively Theorem \ref{ite}) is a straightforward consequence of the combination of \eqref{cramode}, \eqref{crabismode} and \eqref{cromode} (respectively \eqref{oulmode}, \eqref{oulbismode} and \eqref{madomode}) together with the following lemma, which establishes that the residual term $R_n$ (defined as in \eqref{tetemode}) is negligeable.


\begin{lemma}\label{afa}
Let Assumptions (A1)-(A5) hold. If (C2) holds, then $\lim_{n\rightarrow\infty}\tilde{h}_n^{-q}R_n=0$ a.s.
Otherwise, $\lim_{n\rightarrow\infty}\sqrt{n\tilde{h}_n^d}R_n=0$ a.s.
\end{lemma}

\paragraph{Proof of Lemma \ref{afa}}
We first note that a Taylor's expansion implies the existence of $\zeta_n$ such that $\|\zeta_n-\theta_n\|\leq \|\theta_n-\theta\|$ and 
\begin{eqnarray*}
R_n & = & \left(\theta_n-\theta\right)^T\nabla \tilde{f}_n(\zeta_n) \\ & = & \left(\theta_n-\theta\right)^T\left[\nabla \tilde{f}_n(\zeta_n)-\nabla f(\zeta_n) +\nabla f(\zeta_n) -\nabla f(\theta) \right].
\end{eqnarray*}
Let $\mathcal{V}$ be a compact set that contains $\theta$; for $n$ large enough, we get 
\begin{eqnarray*}
\|R_n\| & = & O\left(\|\theta_n-\theta\|\left[\sup_{x \in \mathcal{V}}\|\nabla \tilde{f}_n(x)-\nabla f(x)\|+\|\zeta_n-\theta\|\right]\right) \\ & = & O\left(\|\theta_n-\theta\|\sup_{x \in \mathcal{V}}\|\nabla \tilde{f}_n(x)-\nabla f(x)\|+ \|\theta_n-\theta\|^2\right).
\end{eqnarray*}
On the one hand, let us recall that the a.s. convergence rate of $\left(\theta_n-\theta\right)$ is given by the one of $\left[D^2f(\theta)\right]^{-1}\nabla f_n(\theta)$ (see \eqref{tatamode} and the comment below). One can apply \eqref{oulmode}, \eqref{oulbismode}, and \eqref{madomode} and obtain the exact a.s. convergence rate of $\theta_n-\theta$. However, to avoid assuming (A6), we apply here Lemmas \ref{aea} and \ref{aia} (with $|\alpha|=1$ and $\left(g_n,b_n\right)=(\tilde f_n,\tilde h_n))$, and get the following upper bound of the a.s. convergence rate of $\theta_n-\theta$: for any $\gamma>0$ and $\varepsilon >0$ small enough,
\begin{eqnarray}
\|\theta_n-\theta\| & = & O\left(\sqrt{\frac{(\log n)^{1+\gamma}}{nh_n^{d+2}}}+\frac{\sum_{i=1}^nh_i^q}{n}\right) 
= O\left(\sqrt{\frac{(\log n)^{1+\gamma}}{nh_n^{d+2}}}+h_n^{q-\varepsilon}\right) \quad \mbox{a.s.}
\label{achamode}
\end{eqnarray} 
On the other hand, we have
\begin{eqnarray*}
\sup_{x \in \mathcal{V}}\|\nabla \tilde{f}_n(x)-\nabla f(x)\| & \leq & \sup_{x \in \mathcal{V}}\|\nabla \tilde{f}_n(x)- \mathbb{E}\left(\nabla \tilde{f}_n(x)\right)\|+\sup_{x \in \mathcal{V}}\|\mathbb{E}\left(\nabla \tilde{f}_n(x)\right)-\nabla f(x)\|.
\end{eqnarray*}
The application of Lemmas \ref{aea} and \ref{aia} with $|\alpha|=1$, $\left(g_n,b_n\right)=(\tilde f_n,\tilde h_n)$ ensures that, for any $\gamma>0$ and $\varepsilon>0$ small enough,
\begin{equation}
\sup_{x \in \mathcal{V}}\|\nabla \tilde{f}_n(x)-\nabla f(x)\|  =  O\left(\sqrt{\frac{(\log n)^{1+\gamma}}{n\tilde{h}_n^{d+2}}}+\frac{\sum_{i=1}^n\tilde{h}_i^q}{n}\right) 
=
O\left(\sqrt{\frac{(\log n)^{1+\gamma}}{n\tilde{h}_n^{d+2}}} + \tilde{h}_n^{q-\varepsilon} \right) \quad \mbox{a.s.}
\label{amymode}
\end{equation}
Let $\mathcal{L}$ denotes a generic slowly varying function that may vary from line to line.
\begin{itemize}
\item Let us first assume that (C1) holds. The application of (\ref{achamode}) and (\ref{amymode}) ensures that for any $\varepsilon>0$ small enough,
$${\sqrt{n\tilde h_n^d}\|\theta_n-\theta\|\sup_{x \in \mathcal{V}}\|\nabla\tilde{f}_n(x)-\nabla f(x)\|}  = 
O\left(\mathcal{L}(n)\left[n^{-\frac{1}{2}\left(1-a(d+2)-2\tilde{a}\right)}+n^{\tilde{a}-a(q-\varepsilon)}\right]\right) + o(1) \quad \mbox{a.s.} $$
Observe that by (C1)i), it is straightforward to see that $2\tilde{a}+a(d+2)<1$ and $\tilde{a}<a(q-\varepsilon)$ for any $\varepsilon>0$ small enough, so that 
$\sqrt{n\tilde h_n^d}\|\theta_n-\theta\|\sup_{x \in \mathcal{V}}\|\nabla\tilde{f}_n(x)-\nabla f(x)\|=o(1)$ a.s.
Moreover, the application of (\ref{achamode}) ensures that
\begin{eqnarray*}
\sqrt{n\tilde h_n^d}\|\theta_n-\theta\|^2  & = &
O\left(\mathcal{L}(n)\left[n^{-\frac{1}{2}\left(1-2a(d+2)+\tilde{a}d\right)}+n^{\frac{1}{2}\left(1-\tilde{a}d-4a(q-\varepsilon)\right)}\right]\right)  \quad \mbox{a.s.} 
\end{eqnarray*}
Now, by (C1)ii) we have $2a(d+2)-\tilde{a}d<1$  and $\tilde{a}d+4a(q-\varepsilon)>1$ for any $\varepsilon >0$ small enough, and thus it follows that
$\sqrt{n\tilde h_n^d}\|\theta_n-\theta\|^2 = o(1)$ a.s.,
which ensures the first part of Lemma \ref{afa}.
\item We now assume that (C2) holds. 
Since $\tilde a q\leq {q}/({d+2q})<1$, using \eqref{relmode}, \eqref{achamode} and \eqref{amymode}, we have
\begin{eqnarray}
\label{faumode}
\lefteqn{\frac{1}{\tilde{h}_n^q}\|\theta_n-\theta\|\sup_{x \in \mathcal{V}}\|\nabla\tilde{f}_n(x)-\nabla f(x)\|} \nonumber\\   & = &
O\left(\mathcal{L}(n)\left[n^{-1+\frac{a(d+2)}{2}+\frac{\tilde{a}(d+2q+2)}{2}}+n^{-\frac{1}{2}-a(q-\varepsilon)+\frac{\tilde{a}(d+2q+2)}{2}}\right]\right)+ o(1) \quad \mbox{a.s.} 
\end{eqnarray}
On the one hand, for any $\varepsilon >0$ small enough, it is straightforward to see that condition (C2) implies the following inequalities:
\begin{eqnarray}
\label{ine1mode}
a(d+2)+\tilde{a}(d+2q+2)<2\ \ \mbox{and}\ \  \tilde{a}(d+2+2q)<1+2a(q-\varepsilon),\\
\tilde{a}q+a(d+2)<1 \ \ \mbox{and} \ \ \tilde{a}q<2a(q-\varepsilon).
\label{ine2mode}
\end{eqnarray}
Therefore, it follows from \eqref{faumode} and \eqref{ine1mode} that
\begin{eqnarray*}
\frac{1}{\tilde{h}_n^q}\|\theta_n-\theta\|\sup_{x \in \mathcal{V}}\|\nabla\tilde{f}_n(x)-\nabla f(x)\| & = & o(1)\ \ \mbox{a.s.} 
\end{eqnarray*}
On the other hand, observe again that by (\ref{achamode}) and \eqref{ine2mode}, we have 
\begin{eqnarray*}
\frac{1}{\tilde{h}_n^q}\|\theta_n-\theta\|^2 & = & 
O\left(\mathcal{L}(n)\left[n^{-\left(1-\tilde{a}q-a(d+2)\right)}+n^{\tilde{a}q-2a(q-\varepsilon)}\right]\right) =  
o(1)\quad \mbox{a.s.},
\end{eqnarray*}
which concludes the proof of Lemma \ref{afa}. 
\end{itemize}

\paragraph{Acknowledgements} We deeply thank two anonymous Referees for their 
helpfull suggestions and comments.


\begin{thebibliography}{99}
\bibitem{abraham1} Abraham, C., Biau, G. and Cadre, B. (2003), Simple estimation of the mode of a multivariate density. \textit{Canadian J. Statist.} \textbf{31}, pp 23-34. 
\bibitem{abraham2} Abraham, C., Biau, G. and Cadre, B. (2004), On the asymptotic properties of a simple estimate of the mode. \textit{ESAIM Prob. and Statist.} \textbf{8}, pp 1-11. 
\bibitem{daviesm} Davies, H.I. (1973), Strong consistency of a sequential estimator of a probability density function. \textit{Bull. Math. Statist.}  \textbf{15}, pp. 49-54. 
\bibitem{devroyem} Devroye, L. (1979), On the pointwise and integral convergence of recursive kernel estimates of probability densities. \textit{Utilitas Math.}  \textbf{15}, pp. 113-128.
\bibitem{eddy1m} Eddy, W.F. (1980),
Optimum kernel estimates of the mode,
{\em Ann. Statist.} {\bf 8}, pp. 870-882.

\bibitem{eddy2m} Eddy, W.F. (1982),
The asymptotic distributions of kernel estimators of the mode,
{\em Z. Warsch. Verw. Gebiete} {\bf 59}, pp. 279-290.

\bibitem{einmahlm} Einmahl, U. (1987), A useful estimate in the multidimensional invariance principle. \textit{Probability theory and related fields} \textbf{76}, pp. 81-101. 
\bibitem{fellerm} Feller, W. (1970), An introduction to probability theory and its applications. Second edition Volume II, Wiley.
\bibitem{hall3m} Hall, P. (1992),
Effect of bias estimation on coverage accuracy of bootstrap 
confidence intervals for a probability density,
{\em Ann. Statist.} {\bf 20}, pp. 675-694.
\bibitem{konakovm} Konakov, V.D. (1973),
On asymptotic normality of the sample mode of multivariate distributions,
{\em Theory Probab. Appl.} {\bf 18}, pp. 836-842.
\bibitem{kovalm} Koval, V. (2002), A new law of the iterated logarithm in $\mathbb{R}^d$ with application to matrix-normalized sums of randoms vectors. \textit{Journal of Theoretical Probability}  \textbf{15}, pp. 249-257.  
\bibitem{menonprasadsinghm} Menon, V.V., Prasad, B. and Singh, R.S. (1984), Non-parametric recursive estimates of a probability density function and its derivatives. \textit{Journal of Statistical Planning and inference}  \textbf{9}, pp. 73-82.
\bibitem{mokkadempelletierm} Mokkadem, A. and Pelletier, M. (2003), The law of the iterated logarithm for the multivariate kernel mode estimator. \textit{ESAIM: Probab. Statist.} \textbf{7}, pp. 1-21.
\bibitem{mokka} Mokkadem, A. and  Pelletier, M. (2007),
{A companion for the Kiefer-Wolfowitz-Blum stochastic 
approximation algorithm,} {\em Ann. Statist.,} {\bf 35},  1749-1772.

\bibitem{mokkadempelletierthiamm} Mokkadem, A., Pelletier, M. and Thiam, B. (2006), 
Large and moderate deviations principles for recursive kernel estimators of a multivariate 
density and its partial derivatives. 
\textit{Serdica Math. J.} \textbf{32}, pp. 323-354.

\bibitem{muller} M\"uller H.G. (1989) Adaptive nonparametric peak estimation 
{\em Ann. Statist.,} {\bf 17},  1053-1069.

\bibitem{nadarayam} Nadaraya, E.A. (1965), On non-parametric estimates of density functions and regression curves. \textit{Theory Probab. Appl.} \textbf{10}, pp. 186-190. 
\bibitem{parzenm} Parzen, E. (1962), On estimation of a probability density function and mode. \textit{Ann. Math. Statist.} \textbf{33}, pp. 1065-1076.
\bibitem{petrovm} Petrov, V.V. (1995), Limit theorems in probability theory, Clarendon Press, Oxford.
\bibitem{romanom} Romano, J. (1988), On weak convergence and optimality of kernel density estimates of the mode. \textit{Ann. Statist.} \textbf{16}, pp. 629-647.
\bibitem{rosenblattm} Rosenblatt, M. (1956), Remarks on some-non-parametric estimates of density function. \textit{Ann. Math. Statist.}  \textbf{27}, pp. 832-837. 
\bibitem{roussasm} Roussas, G. (1992), Exact rates of almost sure convergence of a recursive kernel estimate of a probability density function: Application to regression and hazard rate estimate. \textit{J. of Nonparam. Statist.}  \textbf{3}, pp. 171-195.
\bibitem{ruschendorfm} R\"uschendorf, L. (1977), Consistency of estimators for multivariate density functions and for the mode. \textit{Sankhya Ser. A}, \textbf{39}, pp. 243-250.
\bibitem{samantam} Samanta, M. (1973),
Nonparametric estimation of the mode of a multivariate density,
{\em South African Statist. J.}
{\bf 7}, pp. 109-117.

\bibitem{tsybakov} Tsybakov, A.B.  (1990),
Recurrent estimation of the mode of a multidimensional distribution,
{\em Problems Inform. Transmission} {\bf 26}, 31-37

\bibitem{vanryzinm} Van Ryzin, J. (1969), On strong consistency of density estimates. \textit{Ann. Math. Statist.} \textbf{40}, pp 1765-1772.
\bibitem{vieum} Vieu, P. (1996),
A note on density mode estimation,
{\em Statist. Probab. Lett.} {\bf 26}, 297-307

\bibitem{wegmandaviesm} Wegman, E.J. and Davies, H.I. (1979), Remarks on some recursive estimators of a probability density. \textit{Ann. Statist.}  \textbf{7}, pp. 316-327. 
\bibitem{wertzm} Wertz, W. (1985), Sequential and recursive estimators of the probability density. \textit{Statistics}  \textbf{16}, pp. 277-295. 
\bibitem{wolvertonm} Wolverton, C.T. and Wagner, T.J. (1969), Asymptotically optimal discriminant functions for pattern classification. \textit{IEEE Trans. Inform. Theory}  \textbf{15}, pp. 258-265. 
\bibitem{yamatom} Yamato, H. (1971), Sequential estimation of a continuous probability density function and mode. \textit{Bull. Math. Satist.}  \textbf{14}, pp. 1-12.   
\bibitem{ziegler a} Ziegler, K. (2003), On the asymptotic normality of kernel regression estimators 
of the mode in the random design model. 
\textit{J. Statist. Plann. Inf.}  \textbf{115}, pp. 123-144.   
\bibitem{ziegler b} Ziegler, K. (2004), Adaptive kernel estimation of the mode in a nonparametric 
random design regression model.
\textit{Probab. Math. Satist.}  \textbf{24}, pp. 213-235.   
\end{thebibliography}
\end{document}